\nonstopmode
\documentclass[12pt]{amsart}
\usepackage{amscd}
\usepackage{amssymb}
\input xy
\xyoption{all}

\numberwithin{equation}{section}

\swapnumbers

\theoremstyle{plain} 
\newtheorem{thm}[equation]{Theorem}

\newtheorem{lem}[equation]{Lemma}
\newtheorem{prop}[equation]{Proposition}

\newtheorem{axiome}[equation]{Axiom E}
\newtheorem{axiomone}[equation]{Axiom~I}
\newtheorem{axiomtwo}[equation]{Axiom~II}
\newtheorem{axiomthree}[equation]{Axiom~III}
\newtheorem{conditiond}[equation]{Axiom IV}
\newtheorem{conditionf}[equation]{Axiom V}
\newtheorem{axiomsix}[equation]{Axiom VI}

\theoremstyle{definition}
\newtheorem{defn}[equation]{Definition}

\theoremstyle{remark}
\newtheorem{rem}[equation]{Remark}

\newtheorem*{VariableNoNum}{{\VariableText}}
\newtheorem{Variable}[equation]{{\VariableText}}

\theoremstyle{definition}
\newtheorem*{VariableNoNumBold}{{\VariableText}}
\newtheorem{VariableBold}[equation]{{\VariableText}}

\newenvironment{titled}[1]
     {\def\VariableText{{#1}}\begin{VariableNoNum}}
     {\end{VariableNoNum}}

\newenvironment{numbered}[1]
     {\def\VariableText{{#1}}\begin{Variable}}
     {\end{Variable}}

\newenvironment{NumberedSubSection}[1]
     {\def\VariableText{{\textrm{\textbf{#1}}}}\begin{VariableBold}}
     {\end{VariableBold}}

\newlength{\asidelength}

\def\Changed/{\ifvmode\else\vadjust{%
\vbox to 0pt{\vskip -\baselineskip%
\hbox to 0pt{\hss\vrule height 0pt depth 1.2\baselineskip\hskip 1em}\vss}}\fi}
\def\CHanged{\ifvmode\else\vadjust{%
\vbox to 0pt{\vskip -\baselineskip%
\hbox to 0pt{\hss\vrule height 0pt depth 1.2\baselineskip\hskip 1em}\vss}}\fi}

\def\Math#1{\def\MathString{#1}\futurelet\MathDelim\MathChoose}
\def\MathChoose{\ifmmode\let\MathDo\MathString%
              \else\let\MathDo\MathSkip\fi%
              \MathDo}
\def\MathSkip{\ifx\MathDelim/\def\MathDo{$\MathString$\EatOne}%
              \else\def\MathDo{$\MathString$}\fi%
              \MathDo}

\def\Text#1{\def\TextString{#1}\futurelet\TextDelim\TextSkip}
\def\TextSkip{\ifx\TextDelim/\def\TextDo{\TextString\EatOne}%
              \else\let\TextDo\TextString\fi%
              \TextDo}
\def\EatOne#1{}

\def\SkipToEndScan#1\EndScan{}
\def\Scan#1#2#3{\ifx#1#2#3\expandafter\SkipToEndScan\fi\Scan#1}
\def\Upper#1{%
\Scan#1aAbBcCdDeEfFgGhHiIjJkKlLmMnNoOpPqQrRsStTuUvVwWxXyYzZ#1#1\EndScan}
\def\Phrase#1 #2/#3/#4=#5 #6/#7/#8.{%
\expandafter\edef\csname#2#3\endcsname{\noexpand\Text{#6#7}}
\expandafter\edef\csname\Upper#2#3\endcsname{\noexpand\Text{\Upper#6#7}}
\expandafter\edef\csname#1#2#3\endcsname{\noexpand\Text{#5 #6#7}}
\expandafter\edef\csname\Upper#1#2#3\endcsname{\noexpand\Text{\Upper#5 #6#7}}
\expandafter\edef\csname#2#4\endcsname{\noexpand\Text{#6#8}}
\expandafter\edef\csname\Upper#2#4\endcsname{\noexpand\Text{\Upper#6#8}}
}

\newcommand{\Cat}[1]{\mathbf{#1}}
\newcommand{\op}{^{\operatorname{op}}}

\newcommand{\R}{\mathbb{R}}

\newcommand{\RightArrow}[1]{\xrightarrow{#1}} 

\newcommand{\Map}{\operatorname{Map}}

\newcommand{\hhh}{\operatorname{h}\!}
\newcommand{\thhh}{\tilde{\operatorname{h}}}
\def\HomotopyOrbit#1on#2/{\ensuremath{#2_{\hhh#1}}}
\def\RedHomotopyOrbit#1on#2/{\ensuremath{#2_{\thhh#1}}}

\newcommand{\Konof}[1]{\Math{\mathcal{K}_{#1}}}
\newcommand{\Konm}{\Math{\Konof m}}
\newcommand{\Km}{\Math{\Konm}}

\newcommand{\assoc}{\Math{\mathcal{A}}}
\newcommand{\oo}{\Math{\mathcal{O}}}
\newcommand{\cso}{\Math{\oo^{\bullet}}}
\newcommand{\om}{\Math{\omega}}
\newcommand{\holim}{\operatorname{holim}}
\newcommand{\weq}{\sim}
\newcommand{\maph}{\operatorname{Map}^{\operatorname{h}}}

\newcommand{\catOf}[1]{\mathbf{#1}}

\newcommand{\spaces}{\sspace}
\newcommand{\sspace}{\Math{\catOf{Sp}}}
\newcommand{\sspaces}{\sspace}
\newcommand{\nspace}{\Math{\catOf{S}}}

\newcommand{\nspaces}{\nspace}
\newcommand{\Sing}{\operatorname{Sing}}

\newcommand{\myunit}{e}
\newcommand{\gtensor}{\Math{\odot}}
\newcommand{\gunit}{\Math{\myunit_{\gtensor}}}
\newcommand{\ctensor}{\Math{\circ}}
\newcommand{\cunit}{\Math{\myunit_{\textstyle\ctensor}}}
\newcommand{\gentensor}{\Math{\diamond}}
\newcommand{\genunit}{\Math{\myunit_{\gentensor}}}

\newcommand{\operads}{\Math{\catOf{O}}}
\newcommand{\opa}{\Math{\mathcal{P}}}
\newcommand{\opb}{\Math{\mathcal{Q}}}
\newcommand{\cof}{^{\operatorname{c}}}

\newcommand{\fib}{^{\operatorname{f}}}
\Phrase a reduced//=a reduced//.
\newcommand{\Of}[1]{_{#1}}
\newcommand{\monunit}{\epsilon}
\newcommand{\monmult}{\mu}

\newcommand{\origcatM}{\Math{\catOf{C}}}
\newcommand{\catM}{\origcatM}
\newcommand{\catMon}{\Math{\catOf{M}}}

\newcommand{\bimod}[2]{\text{$#1{-}#2$-bimodule}}

\newcommand{\UnPointedCat}[1]{#1}
\newcommand{\DistSign}{{\operatorname{d}}}
\newcommand{\MaybeDistSign}{{\operatorname{pd}}}
\newcommand{\PointedCat}[1]{#1^{\myunit}}
\newcommand{\ModuleCat}[2]{#1_{#2}}

\newcommand{\BiModuleCat}[3]{#1_{#2,#3}}
\newcommand{\DistModuleCat}[2]{#1_{#2}^\DistSign}
\newcommand{\DistBiModuleCat}[3]%
        {#1_{#2,#3}^{\DistSign}}
\newcommand{\MaybeDistBiModuleCat}[3]%
        {#1_{#2,#3}^{\MaybeDistSign}}

\newcommand{\PointedModuleCat}[2]{\ModuleCat{\PointedCat{#1}}{#2}}
\newcommand{\PointedBiModuleCat}[3]%
        {\BiModuleCat{\PointedCat{#1}}{#2}{#3}}
\newcommand{\PointedDistModuleCat}[2]%
        {\DistModuleCat{\UnPointedCat{#1}}{#2}}
\newcommand{\PointedDistBiModuleCat}[3]%
        {\DistBiModuleCat{\UnPointedCat{#1}}{#2}{#3}}
\newcommand{\MaybeDistModuleCat}[2]{#1_{#2}^\MaybeDistSign}
\newcommand{\PointedMaybeDistBimoduleCat}[3]%
        {\MaybeDistBiModuleCat{\UnPointedCat{#1}}{#2}{#3}}  

\newcommand{\catBirs}{\Math{\BiModuleCat{\catM}{R}{S}}}
\newcommand{\catModrs}{\catBirs}
\newcommand{\catBirr}{\Math{\BiModuleCat{\catM}{R}{R}}}
\newcommand{\catBirsp}{\Math{\PointedBiModuleCat{\catM}{R}{S}}}
\newcommand{\catBirspd}{\Math{\PointedDistBiModuleCat{\catM}{R}{S}}}
\newcommand{\catBirsd}{\Math{\MaybeDistBiModuleCat{\catM}{R}{S}}}
\newcommand{\catMonrs}{\Math{\BiModuleCat{{\catMon}}{R}{S}}}
\newcommand{\catMods}{\Math{\ModuleCat{\catM}{S}}}
\newcommand{\catModsp}{\Math{\PointedModuleCat{\catM}{S}}}
\newcommand{\catModsd}{\Math{\MaybeDistModuleCat{\catM}{S}}}

\newcommand{\catBiab}{\Math{\BiModuleCat{\catM}{\opa}{\opb}}}

\newcommand{\catModab}{\catBiab}

\newcommand{\catBiabp}{\Math{\PointedBiModuleCat{\catM}{\opa}{\opb}}}
\newcommand{\catBiaap}{\Math{\PointedBiModuleCat{\catM}{\opa}{\opa}}}

\newcommand{\catMonab}{\Math{\BiModuleCat{{\catMon}}{\opa}{\opb}}}
\newcommand{\catMonaa}{\Math{\BiModuleCat{{\catMon}}{\opa}{\opa}}}
\newcommand{\catModb}{\Math{\ModuleCat{\catM}{\opb}}}

\newcommand{\catD}{\Math{\catOf D}}
\newcommand{\catC}{\Math{\catOf C}}

\newcommand{\catMonrsd}{\Math{\DistBiModuleCat{\catMon}{R}{S}}}

\newcommand{\Envel}{E}
\newcommand{\Envelh}{\Envel^{\operatorname{h}}}

\newcommand{\forget}{\natural}

\newcommand{\mapbih}[2]{\maph_{#1{-}#2}}
\newcommand{\maptrih}[3]{\maph_{#1{-}#2,#3}}

\newcommand{\gassoc}{\Math{A}}

\newcommand{\RightModuleOver}[2]{\Math{#1_{#2}}}
\newcommand{\nspacesRightModule}[1]{\RightModuleOver{\nspaces}{#1}}
\newcommand{\nspacesA}{\RightModuleOver{\nspaces}{\assoc}}

\newcommand{\End}{\operatorname{End}}
\newcommand{\Auth}{\operatorname{Aut}^{\operatorname{h}}}
\newcommand{\auth}{\Auth}

\newcommand{\PotentialDistSign}{\MaybeDistSign}
\newcommand{\Moduli}{\mathcal{M}}
\newcommand{\DistinguishedModuli}[2]{\Math{\Moduli^{\DistSign}(#1,#2)}}
\newcommand{\UndistinguishedModuli}[2]{\Math{\Moduli^{\PotentialDistSign}(#1,#2)}}
\newcommand{\RightModuli}[1]{\Math{\Moduli^{\PotentialDistSign}({\genunit},#1)}}

\newcommand{\mycoprod}{\sqcup}

\Phrase an acceptable//=an acceptable//.

\newcommand{\catMC}{\Math{\catOf{D}}}
\newcommand{\catMCp}{\Math{\catOf{E}}}

\newcommand{\eqto}{\RightArrow{\sim}}

\newcommand{\Overcat}[2]{#2{\scriptstyle\searrow}#1}
\newcommand{\Undercat}[2]{#1{\scriptstyle\searrow}#2}

\newcommand{\rqp}{\rho}
\newcommand{\lqp}{\lambda}

\def\SquigThree{\Math{\catOf{Z}}}

\newcommand{\monmap}{\om}

\Phrase a gradedspace//s=a graded space//s.
\Phrase a gradedsubspace//s=a graded subspace//s.
\Phrase a gradedset//s=a graded set//s.

\newcommand{\underPair}[2]{#1{-}#2}

\newcommand{\freeop}{\Phi}
\newcommand{\pfreeop}{\PointedCat{\Phi}}
\newcommand{\freeopp}{\pfreeop}
\newcommand{\catMp}{\Math{\PointedCat{\catM}}}

\newcommand{\Identity}[1]{1_{#1}}

\newcommand{\height}{\chi}
\newcommand{\lc}{\lambda}
\newcommand{\bc}{\beta}

\newcommand{\tree}{\Math{t}}
\newcommand{\etree}{\Math{e}}
\def\treeof(#1;#2){[#1;#2]}
\newcommand{\Trees}{\mathbb{T}}
\newcommand{\treelab}[2]{#2(#1)}
\newcommand{\Treelab}[1]{\Trees(#1)}
\newcommand{\Tof}{\mathbb{T}}

\Phrase a branchcount//s=a branching count//s.

\newcommand{\comp}{F}
\newcommand{\eval}{\epsilon}

\newcommand{\ReflexArrow}{\Rightarrow}
\newcommand{\seqcoprod}{\cup}
\newcommand{\gradedspacesum}{\seqcoprod}
\newcommand{\genof}[1]{[#1]}
\newcommand{\shape}[1]{\operatorname{sh}(#1)}
\newcommand{\labset}[1]{\operatorname{Lab}(#1)}

\newcommand{\opaa}{\Math{\mathcal{R}}}
\newcommand{\opbb}{\Math{\mathcal{S}}}
\newcommand{\Treelabpair}[3]{\Tof_{#3}(#1,#2)}

\newcommand{\colltreelab}[3]{\Tof_{#3}^{\operatorname{c}}(#1,#2)}
\newcommand{\collTreelab}{\colltreelab}

\newcommand{\Filt}[1]{\operatorname{Filt}_{#1}}

\newcommand{\colim}{\operatorname{colim}}

\newcommand{\longX}{\Math{X_1^+}}
\newcommand{\shortX}{\Math{X_1^-}}
\newcommand{\bothX}{\Math{X_1^\pm}}
\newcommand{\bothd}{\Math{d^\pm}}

\newcommand{\simsub}[1]{\sim_{#1}}
\newcommand{\hres}{\operatorname{H}}
\newcommand{\dhres}{\operatorname{dH}}

\newcommand{\dehres}{d\Envel\!\hres}

\newcommand{\circpower}[2]{#1^{\ctensor #2}}
\newcommand{\diag}{\operatorname{diag}}
\newcommand{\aug}{\alpha}
\newcommand{\daug}{\aug}

\newcommand{\simptensor}{\cdot}
\newcommand{\triple}[3]{(#1,#2,#3)}
\newcommand{\peh}{\Envelh}

\newcommand{\eend}{\Math{\mathcal{E}}}
\newcommand{\eendb}{\Math{\eend_{\opb}}}

\Phrase a classicmonoid//s=a classic monoid//s.
\Phrase a gradedclassicmonoid//s=a graded classic monoid//s.
\Phrase a classicgradedmonoid//s=a graded classic monoid//s.

\newcommand{\gcm}{\Math{\Cat{M}}}
\newcommand{\cgm}{\gcm}
\newcommand{\Maph}{\maph}

\newcommand{\noRtA}{^{\prime}}
\newcommand{\gassocn}{\Math{\gassoc\noRtA}}

\newcommand{\RtAmod}{\Math{\Cat{T}}}

\newcommand{\LeftAmod}{\Math{\Cat{M}}}
\newcommand{\LeftAAmod}{\Math{\Cat N}}
\newcommand{\FunnyBiMods}{\Math{\BiModuleCat{\RtAmod}{\gassoc}\gassoc}}

\newcommand{\coprodOver}[1]{\mycoprod^{#1}}
\newcommand{\gcma}{\Math{G}}
\newcommand{\gcmb}{\Math{H}}

\newcommand{\cd}{\Math{\Cat D}}

\newcommand{\hocolim}{\operatorname{hocolim}}
\begin{document}
\title[Long knots]{Long knots and maps between operads}

\author{William Dwyer}
\author{Kathryn Hess}

\address{Department of Mathematics, University of Notre Dame, Notre
  Dame, IN 46556, USA}
\email{dwyer.1@nd.edu}

\address{Institut de G\'eom\'etrie, alg\`ebre et topologie (IGAT),
  \'Ecole Polytechnique F\'ed\'erale de Lausanne, CH-1015 Lausanne, Switzerland}

\email{kathryn.hess@epfl.ch}

\date{\today}
\thanks{WGD was partially supported by NSF grant DMS~0354787}
\begin{abstract}
  We identify the space of tangentially straightened long knots in $\mathbb{R}^m$,
  $m\ge4$, as the double loops  on the space of derived operad maps from
  the associative operad into a version of the little $m$-disk operad.
  This verifies a conjecture of Kontsevich, Lambrechts, and Turchin.
\end{abstract}
\maketitle

\section{Introduction}\label{CIntroduction}
A \emph{long knot} in Euclidean $m$-space $\R^m$ is
a smooth embedding $\R\to\R^m$ which agrees with inclusion of the
first coordinate axis on the complement of some compact set in~$\R$; a
\emph{tangential straightening} for the knot is a null homotopy (constant
near~$\infty$) of the map
$\mathbb{R}\to S^{m-1}$ obtained by taking the unit tangent vector of
the knot at each point. See \cite[5.1]{sinha} for more details.
Starting from work of Goodwillie, Klein, and Weiss \cite{weissPoint}
\cite{weissGoodwilliePoint} \cite{goodwillieEstimates}
\cite{goodwillieKleinMultiple} on embedding spaces, Sinha \cite{sinha} has
proved that the space of tangentially straightened long knots in $\R^m$ is
equivalent to the homotopy limit of an explicit cosimplicial space
constructed from the $m$'th Kontsevich operad~\Konm. 

Let \assoc/ denote the associative operad. The cosimplicial space
Sinha considers is derived by formulas of McClure and Smith
\cite{mcclureSmith} from an operad map $\assoc\to\Km$; more generally,
McClure and Smith build a cosimplicial space \cso/ from {any}
operad \oo/ of spaces and operad map $\assoc\to\oo$. Our main theorem 
gives an independent formula for the homotopy limit of this construction.  Say that the operad \oo/
is \emph{\reduced/} if $\oo\Of0$ and $\oo\Of1$ are weakly contractible.

\begin{thm}\label{MainTheorem} (\S\ref{CFinalProof})
  Let \oo/ be \areduced/ operad of simplicial sets, $\om:\assoc\to\oo$ an
  operad map, and \cso/ the cosimplicial object associated to $\om$ by
  \cite{mcclureSmith}. Then there is a natural weak homotopy
  equivalence
  \[
      \holim \cso \weq \Omega^2\maph(\assoc,\oo)_{\om}\,.
  \]
\end{thm}

The space on the right is the double loops, based at~\om, on a
derived space of operad maps $\assoc\to\oo$. See below for more
details.  Let $\Sing(-)$ be
the singular complex functor. Since \Km/ is reduced, applying
\ref{MainTheorem} with $\oo=\Sing(\Km)$ shows
that the space of tangentially straightened long knots in $\R^m$ is weakly homotopy
equivalent to the double loop space on a space of operad maps from the
associative operad to $\Sing(\Km)$. This verifies a conjecture of
Kontsevich \cite[2.17{\it ff\/}]{sinha} as adjusted by Lambrechts and
Turchin and re-expressed simplicially.  The operad \Km/ is
a pared-down model of the little $m$-disk operad
\cite{kontsevichOperads} \cite[\S4]{sinha}. A result similar to
\ref{MainTheorem} has been announced recently by Lambrechts and Turchin.

\begin{NumberedSubSection}{More details}\label{MoreDetails}
  From this point on in the paper ``space'' means ``simplicial set.''
  Let \sspace/ denote the category of spaces and \nspace/
  the category of \gradedspaces/ (sequences of spaces indexed by the nonnegative
  integers).
  The \emph{composition product} $X\ctensor Y$ of \gradedspaces/  is
  given by
  \begin{equation}\label{DefineCompositionProduct}
             (X\circ Y)\Of n = \coprod_{i,\,\,j_1+\cdots j_i=n} X\Of i\times
             Y\Of{j_1}\times\cdots\times Y\Of{j_i}\,.
  \end{equation}
  An \emph{operad} is a monoid object for this associative (but nonsymmetric) product.
   See \cite{mayGeometry} or \cite[2.13]{sinha}, but
  note that we do not impose constraints on an operad in levels~0
  and~1.   The \emph{associative operad} $\assoc$ is
  the \gradedspace/  which has a single point at every level $n$,
  $n\ge0$; it has a unique  operad structure. We will let \operads/
  denote the category of operads.

  The category \nspace/ is a simplicial model category in
  which a map $X\to Y$ is an equivalence, cofibration, or fibration
  if and only if each individual map $X\Of n\to Y\Of n$ has the corresponding
  property in the usual model structure on spaces \cite[3.2]{hovey}. It follows from \cite{rezkThesis}
  (see also \cite[pf. of 7.2]{rezkProperModel}) that \operads/ is a
  simplicial model category in which a map $\opa\to\opb$ is
  an equivalence or a fibration if and only if the underlying map in
  \nspaces/ has the corresponding property. The derived mapping space
  $\maph(\assoc,\oo)$ in \ref{MainTheorem} is the
  simplicial mapping space $\Map_{\operads}(\assoc\cof,\oo\fib)$,
  where $\assoc\cof$ is a cofibrant replacement for \assoc/ and
  $\oo\fib$ a fibrant replacement for \oo.

  In describing the left hand side of \ref{MainTheorem}, it is
  convenient to write a point $f\in\oo\Of i$ as if it were a function
  $f(x_1,\cdots, x_i)$ of $i$~formal variables. If the  image $h$ of a tuple
  $(f,g_1,\cdots, g_i)$  under the operad structure map
  $\oo\ctensor\oo\to\oo$ is written as a composite
 \[
    h(x_1,\cdots,x_{J}) =
    f(g_1(x_1,\cdots,x_{j_1}),\ldots,
    g_i(x_{J-j_i+1},\cdots,x_{J}))
  \]
  (where $g_k\in\oo\Of{j_k}$ and $J=\sum j_k$), then the operad
  identities for the structure map express the associativity of
  composition.  (Writing~$h$ as a composite in this way amounts to
  pretending for the sake of notation that $\oo$ is the endomorphism
  operad of a space.)  Given a map $\assoc\to\oo$, let
  $\monunit\in\oo\Of0$ and $\monmult\in\oo\Of2$ be the images
  respectively of $\assoc\Of0\to\oo\Of0$ and $\assoc\Of2\to\oo\Of2$.
  We will think of $\monmult$ as a multiplication and temporarily
  write $x*y$ for $\monmult(x,y)$.  For each $n\ge0$ the cosimplicial
  space $\cso$ has $\oo\Of n$ in cosimplicial degree~$n$; the coface
  and codegeneracy operators on $f\in\oo\Of n$ are given by
  \[
      (d^i f)(x_1,\ldots,x_{n+1}) = 
           \begin{cases}
                 x_1*f(x_2,\ldots,x_{n+1}) & i=0\\  
                 f(x_i,\ldots,x_i*x_{i+1},\ldots x_{n+1}) & 1\le i\le
                 n\\
                 f(x_1,\ldots,x_n)*x_{n+1} & i=n+1
             \end{cases}
   \]
  \[
       (s_if)(x_1,\ldots x_{n-1}) =
       f(x_1,\ldots,x_{i-1},\monunit,x_{i},\ldots,x_{n-1})
  \]
  For other descriptions of \cso/ see \cite[2.17]{sinha}
  or \cite[\S3]{mcclureSmith}. The above formulas
  recall a Hochschild resolution or
  a two-sided bar construction, and this turns out not to be
  accidental (see the proof of \ref{IdentifyCoSimplicial}).

\end{NumberedSubSection}

\begin{NumberedSubSection}{Method of proof} Oddly enough, our proof of
  \ref{MainTheorem} comes down to applying a single principle twice,
  in very different cases; each application gives rise to
  one instance of $\Omega$. Before looking at these applications, we will
  describe the principle.
  
  \begin{numbered}{A connection between maps and bimodules}\label{Connection}
  Suppose that $(\catM,\gentensor,\genunit)$ is a category with a monoidal product
  $\gentensor$ for which $\genunit$ is the unit, and that $\catMon$ is the category of monoids in
  \catM. It is not
  necessary for $\gentensor$ to be symmetric
  monoidal. 
Given monoids $R$, $S$ in \catM, the notions of left $R$-module, right $S$-module, and \bimod
  RS are defined as usual. A
  \emph{pointed} module or bimodule $X$ is one which is 
  provided with a \hbox{\catM-map} $\genunit\to X$; a \emph{monoid under}
  $\underPair RS$ is a monoid $T$ together with monoid maps $R\to T$ and $S\to T$.

  Suppose that there is
  a model category structure on \catM/ which induces compatible (\S\ref{MonoidalCategories})
  model structures on all monoid and bimodule categories.
  There is a  forgetful functor
  $\forget$ from monoids under $\underPair RS$ to pointed \bimod RSs, with
  the basepoint for $\forget T$ given by the unit map $\genunit\to T$.
  Assume that $\forget$ has a left
  adjoint $\Envel$ (the \emph{enveloping monoid} functor) with the property that $(\Envel,\forget)$ forms a
  Quillen pair, and  
 let $\Envelh$ denote the left derived functor~(\ref{AutoQuillenPair})
 of~$\Envel$~. 

Say that a pointed right $S$-module~$X$ is 
\emph{distinguished}  if the map 
\[S\cong
  \genunit\gentensor S\to X\gentensor S\to X
\]
induced by $\genunit\to X$ and $X\gentensor S\to S$ is an equivalence;
an \bimod RS is {distinguished} if it is distinguished as
a right $S$-module. Similarly a monoid $T$ under $\underPair RS$
is distinguished if $T$ is distinguished as a right $S$-module, or
equivalently, if the structure map $S\to T$ is an equivalence. 
There is one major axiom:

\begin{axiome}\label{AxiomE}
  If $X$ is a
  distinguished  \bimod RS, then $\Envelh(X)$ is a
  distinguished monoid under~$\underPair RS$.
\end{axiome}

The next theorem is the source of~$\Omega$.

\begin{thm}\label{GeneralTheorem} (\ref{LoopedMonoidMapFibration})
  Let $\om:R\to S$ be a map of monoids in \catM.
Then under the above assumptions, and some additional technical
  conditions, there is a fibration sequence
\[
  \Omega \maph_{\catMon}(R,S)_{\om} \to \mapbih
  RR(R,S)\to\maph_{\catM}(\genunit, S)\,.
 \]
\end{thm}

In this statement $\maph$ stands for the general model category
theoretic derived mapping space (\ref{ModelGeneralities}), which agrees with the appropriate derived
simplicial mapping space if the model category involved is a
simplicial model category. The middle space is computed in the
category of \bimod RRs, with $S$ treated as an \bimod RR
via~\om. The right hand map is induced by the unit map $\genunit\to
R$, and the homotopy fibre is meant to be computed over the unit map
$\genunit\to S$. 
\end{numbered}
  
\begin{numbered}{The first application of \ref{GeneralTheorem}} Here
  $(\catM,\gentensor,\genunit)$ is $(\nspaces,\ctensor,\cunit)$; the
  unit $\cunit$ is 
  a \gradedspace/ which is  empty  except for a
  single point at level~1. The monoids in \catM/  are the operads. For any \gradedspace/ $X$, the mapping space
  $\maph_{\nspaces}(\cunit,X)$ is equivalent to $X\Of 1$, and so
  \ref{GeneralTheorem} gives
  the following.

  \begin{thm}\label{ApplicationOne} (\S\ref{CMapsBetweenOperads})
   Suppose that $\om:\assoc\to \oo$ is a map of operads,
    with $\oo\Of1$ contractible. Then there is an equivalence
   \[
           \Omega\maph_{\operads}(\assoc,\oo)_{\om}\weq \mapbih
           \assoc\assoc(\assoc, \oo)\,,
    \]
    where on the right $\oo$ is treated as a \bimod \assoc\assoc{} via~$\om$.
  \end{thm}
  
\end{numbered}

\begin{numbered}{The second application of \ref{GeneralTheorem}}\label{GradedTensorProduct} This is more
  peculiar. For any two objects $X$, $Y$ of \nspaces, the \emph{graded
    cartesian product} $X\gtensor Y$ is the \gradedspace/  defined by
\[ (X\gtensor Y)\Of n = \coprod_{i+j=n} X\Of i\times Y\Of j\,.\]
The unit for $\gtensor$ is the \gradedspace/  \gunit/ which is empty except
for a single point at level~0; under the pairing $\gtensor$ the
category $\nspaces$ becomes a \emph{symmetric} monoidal category. 

It is easy to see that \agradedspace/ $X$ is a  $\gtensor$-monoid
 ($X\gtensor X\to X$) if and only if it is a left
\assoc-module ($\assoc\ctensor X\to X$). But more is true. For any
three \gradedspaces/ $X$, $Y$, $Z$ there is a natural distributive isomorphism
\begin{equation}\label{DistributiveLaw}
           (X\gtensor Y)\ctensor Z\cong (X\ctensor Z)\gtensor
           (Y\ctensor Z)\,.
\end{equation}
This guarantees that if $X$ and $Y$ are right modules over an operad
\opa, then $X\gtensor Y$ is also naturally a right module over $\opa$;
indeed, \gtensor/ provides a symmetric monoidal structure on the
category $\nspacesRightModule\opa$ of right \opa-modules, such that the monoids
 in $(\nspacesRightModule\opa,\gtensor,\gunit)$ are exactly the \bimod
\assoc\opa{}s in $(\nspaces,\ctensor,\cunit)$.

In our second application of \ref{GeneralTheorem}, we take
$\opa=\assoc$ in the above remarks and we let
$(\catM,\gentensor,\genunit)$ be $(\nspacesA,\gtensor,\gunit)$. Let
$\gassoc$ denote the associative operad $\assoc$ with an emphasis on  its role as a
monoid in \nspacesA.  Theorem \ref{GeneralTheorem} then translates to this.

\begin{thm}\label{ApplicationTwo} (\S\ref{CProofOfApplicationTwo})
  Suppose that $\alpha:\assoc\to X$ is a map of \bimod\assoc\assoc{}s,
  and  that
  $X\Of0\sim *$. Then there is an equivalence
  \[      \Omega \mapbih\assoc\assoc(\assoc,X)_\alpha\sim
  \maptrih\gassoc\gassoc\assoc(\gassoc, X)\,.
   \]
\end{thm}
   The mapping space on the right is computed in the category of
   \bimod\gassoc\gassoc{}s with respect to \gtensor{} in a setting in
   which all of the \gradedspaces/ involved are right \assoc-modules
   with respect to~\ctensor! Fortunately, this is a lot less
   complicated than it looks. It turns out that there is a very small
   Hochschild resolution of \gassoc/ as an \bimod \gassoc\gassoc{} in
   the category of right $\assoc$-modules (\S\ref{CFinalProof}),  a resolution which
   in conjunction with \ref{ApplicationTwo} and
   \ref{ApplicationOne} leads directly to \ref{MainTheorem}.
\end{numbered}
\end{NumberedSubSection}

\begin{NumberedSubSection}{Some additional comments}  
 Axiom~\ref{AxiomE} is a disguised form of another assumption.  In
  fortunate cases, enhancing a distinguished right $S$-module $X_S$ to an
  \bimod RS amounts to giving a map from $R$ to an enriched
  endomorphism object $\End^+(X_S)$, but
\ref{AxiomE} suggests that such an enhancement amounts to a map $R\to
S$ of monoidal objects. In spirit, then, Axiom~\ref{AxiomE} is the
assumption that $S\sim\End^+(X_S)$ or more explicitly,
given that $X_S$ is distinguished, that $S$ is equivalent to the
enriched endomorphism object of $S$ itself as a right $S$-module. One
feature of \ref{AxiomE} is that it avoids any direct consideration of what
such an enriched endomorphism object might be.

In practice we prove a delooped version of \ref{GeneralTheorem}. Let
$\DistinguishedModuli RS$ denote the \emph{moduli space} of
distinguished \bimod RSs; this is the nerve of the category whose
objects are distinguished \bimod RSs and whose morphisms are the
equivalences between them. This can be identified
(\ref{ModelGeneralities}) as 
\[\DistinguishedModuli RS \sim \coprod_{\{X\}}B\Auth_{R-S}(X) \,,\]
where $\{X\}$ runs over equivalence classes of distinguished \bimod
RSs, and $\Auth_{R-S}(X)$ is the space of derived self-equivalences of~$X$.

\begin{thm}\label{DeloopGeneralTheorem} (cf. \ref{IdentifyMonoidMaps})
  In the situation of  \ref{GeneralTheorem} there is an
  equivalence
  \[   \maph_{\catMon}(R,S) \sim \DistinguishedModuli RS\,.\]
\end{thm}

There is another way to express~\ref{DeloopGeneralTheorem}. Let
$\UndistinguishedModuli RS$ denote the moduli space of all
potentially distinguished \bimod RSs, i.e., bimodules
which are (abstractly) equivalent to $S$ as right $S$-modules, so that
in particular
$\RightModuli S$ is the moduli space of all potentially distinguished
right $S$-modules. This last space is weakly equivalent to
$B\Auth_S(S)$, where the notation indicates  that $\Auth$ is computed in the category of
right $S$-modules.
Unpacking the basepoint from
\ref{DeloopGeneralTheorem} (cf. \ref{MonoidMapFibration}) reveals that there is a fibration
sequence
\begin{equation}\label{VarDeloopMainTheorem}
\maph_{\catMon}(R,S)\to \UndistinguishedModuli RS\to \RightModuli S\,.
\end{equation}
Here are some contexts
$(\catM,\gentensor,\genunit)$ in which \ref{VarDeloopMainTheorem}
comes up. 

\begin{titled}{Simplicial monoids} Here $\catM$ is the category
  \sspaces/ and \gentensor/ is cartesian product. A monoidal object is a
  simplicial monoid. If $G$ and $H$ are two simplicial monoids which
  happen to be simplicial groups, then \ref{VarDeloopMainTheorem} can
  be identified with the fibration sequence
 \[ \Map_*(BG,BH) \to \Map(BG,BH) \to BH\,.\]
\end{titled}

\begin{titled}{Operads} Here $\catM$ is the category of symmetric
  sequences in \sspaces, \gentensor/ is the appropriate analog of
  the composition product, and the monoidal objects are 
  $\Sigma$-operads.  Let \opa/ and \opb/ be two 
  $\Sigma$-operads. Under the additional assumption that \opb/ is an
  endomorphism operad, the fibration \ref{VarDeloopMainTheorem} for
  $\maph_{\catMon}(\opa,\opb)$ appears as \cite[1.1.5]{rezkThesis}.
\end{titled}

\begin{titled}{Ring spectra} Here $\catM$ is the category of spectra,
  $\gentensor$ is smash product, and the monoidal objects are ring
  spectra. Let $R$ and $S$ be two ring spectra.  A calculation of
  $\Omega\maph_{\catMon}(R,S)$ very similar to \ref{GeneralTheorem} appears in
  \cite{lazarev}. In this case \ref{VarDeloopMainTheorem} gives a
  fibration sequence
  \[
          \maph_{\catMon}(R,S)\to \UndistinguishedModuli RS \to BS^\times\,,
  \]
  where $S^\times$ is the group-like simplicial monoid of units in~$S$.
\end{titled}

\end{NumberedSubSection}

\begin{NumberedSubSection}{Notation and
    terminology}\label{Notation}\label{QuillenTheoremB}
  The word \emph{equivalence} usually refers to equivalence in an
  ambient model category; to avoid certain ambiguities, we sometimes
  use \emph{weak equivalence} or \emph{weak homotopy equivalence} to
  refer to an equivalence in the usual model category of simplicial
  sets. We sometimes elide the distinction between a category and its
  nerve, so that a functor is described as a weak equivalence
  if it induces a weak equivalence on nerves. Adjoint
  functors are always weak  equivalences; more generally, a
  functor $F:\catOf C\to\catD$ is a weak equivalence 
if there is a functor $G:\catD\to\catOf
  C$ such that the composites $FG$ and $GF$ are connected to the
  identity functor by zigzags of natural transformations. If $\catOf
  C$ is a category with some notion of equivalence, the \emph{moduli
    category} $\Moduli(\catOf C)$ is the category of equivalences
  in~$\catOf C$. The \emph{moduli space} of~$\catOf C$ (denoted
  identically) is the nerve of the moduli category.
  \label{DefineModuli}

  If $F:\catOf C\to \catD$ is a functor and~$d$ is an object
  of~$\catD$, $\Overcat{d}{F}$ denotes the over category (comma
  category) of $F$ with respect to~$d$. This is the category whose
  objects are pairs $(c,g)$ where $c\in\catOf C$ and $g$ is a map
  $F(c)\to d$ in \catD; a morphism $(c,g)\to (c',g')$ is a map $c\to
  c'$ in $\catC$ rendering the appropriate diagram commutative. The
  corresponding under category is denoted $\Undercat{d}{F}$. If $F$ is
  the identity functor on~\catD, these categories are
  denoted~$\Overcat{d}{\catD}$ and~$\Undercat{d}{\catD}$. 

Suppose that
  $F:\catC\to \catD$ is a functor 
  such that for every $h:d\to d'$ in~\catD/ the map
  $\Overcat{d}{F}\to\Overcat{d'}{F}$ induced by composition with~$h$
  is a weak homotopy equivalence; in this circumstance Quillen's
  Theorem~B guarantees that for any $d\in \catD$ the homotopy fiber of
  (the nerve of) $F$ over the vertex of~\catD/ represented by~$d$ is
  naturally weakly homotopy equivalent to $\Overcat{d}{F}$. A similar
  result holds with over categories replaced by under categories.

    Our symbol for coproduct is usually~$\mycoprod$, or
    $\mycoprod^{\catD}$ if the ambient category is specified;  
    in \S\ref{CModelFreeCoprod} and \S\ref{CMapsBetweenOperads}
    the symbol $\seqcoprod$ is used for the coproduct of \gradedspaces. We refer
    to the \emph{dimension} of a simplex in a space, the \emph{level}
    or \emph{grade}
    of a constituent of \agradedspace, and the simplicial
    \emph{degree} of a constituent of a (co)simplicial (graded) space.
  \end{NumberedSubSection}

  \begin{rem} In
    this paper we work only with non-$\Sigma$ operads, also called
    planar operads. Many of our results apply to $\Sigma$-operads
    (symmetric operads) or
    even  to multicategories, but we decided to leave this generality
    for later. 
  \end{rem}

\begin{titled}{Organization of the paper} 
Section \ref{CModelCategories} sets up some model category machinery,
which is used in \S\ref{MonoidalCategories} to give a proof of
\ref{GeneralTheorem}. Section \ref{CModelFreeCoprod} moves to a more
concrete consideration of model structures on operads of spaces, and
proves the crucial result the coproduct functor preserves
equivalences. Section \ref{CResolutions} describes the Hochschild
resolution of an operad \opa/ as a bimodule over itself, and shows that
applying an enveloping construction to the Hochschild resolution leads to a resolution of \opa/ as
an operad; \S\ref{CDistinguished} uses this last resolution to
give a proof of \ref{ApplicationOne}. The final two sections transpose
the earlier results to the context of classic monoids in the category
of graded spaces, and go
on to deduce \ref{ApplicationTwo} and \ref{MainTheorem}.

For our approach to the homotopy theory of operads we are deeply indebted
to the results on Rezk \cite{rezkProperModel}
    
\end{titled}

\section{Model categories}\label{CModelCategories}
In this section we develop  basic properties of model categories which
we will need later on; \cite{hovey} and
\cite{goerssJardine} are two background references. The underlying 
definition of model category is from
\cite[1.1.4]{hovey}; in particular, a model category is
\emph{closed} in the sense of Quillen \cite{quillen}, has functorial
factorizations \cite[1.1.1]{hovey}, and has all small limits and
colimits.

\begin{NumberedSubSection}{Some generalities}\label{ModelGeneralities}
  Suppose that $\catMC$ is a model category. 
The symbol $\sim$ marks equivalences in \catMC. A (co)fibration is
  \emph{acyclic} if it is also an equivalence. A \emph{cofibrant
    replacement} for an object $X$ is an equivalence $X\cof\eqto X$
  with cofibrant domain, usually
  obtained (functorially) by factoring the map from the initial object
  to $X$ as a cofibration followed by an acyclic fibration. A fibrant
  replacement $X\eqto X\fib$ is constructed similarly. 

\begin{titled}{Quillen pairs and Quillen equivalences}
An adjoint pair $\lqp\colon\catMC\leftrightarrow\catMCp\colon\rqp$ of functors between model categories is
a \emph{Quillen pair} \cite[1.3.1]{hovey} if $\lqp$ preserves cofibrations and
acyclic cofibrations (equivalently, $\rqp$ preserves fibrations and
acyclic fibrations). In this case, $\lqp$ preserves equivalences
between cofibrant objects and $\rqp$ preserves equivalences between
fibrant objects. The pair $(\lqp,\rqp)$ forms a \emph{Quillen
  equivalence} \cite[1.3.3]{hovey}  if
for all cofibrant $X\in\catMC$ and fibrant $Y\in\catMCp$, a map $\lqp
X\to Y$ is an equivanece in $\catMCp$ if and only if the adjoint
$X\to\rqp Y$ is an equivalence in $\catMC$. Let $\catMC\cof$,
$\catMC\fib$, etc., denote appropriate full subcategories of cofibrant
or fibrant objects. Given a Quillen pair $(\lqp,\rqp)$ there is an
induced diagram of moduli spaces
\begin{equation}\label{QPModuliSpaceDiagram}
     \xymatrix{
           {\Moduli(\catMC\cof)} \ar[r]^\sim\ar[d]^{\lqp}  
                 & {\Moduli(\catMC)} &
                       {\Moduli(\catMC\fib)}\ar[l]_\sim \\
          {\Moduli(\catMCp\cof)}  \ar[r] ^\sim& 
                     {\Moduli(\catMCp)} &
                          {\Moduli(\catMCp\fib)}\ar[u]^{\rqp}\ar[l]_\sim
      }
\end{equation}
in which (by functorial fibrant/cofibrant replacement constructions) the indicated arrows
are homotopy equivalences with explicit homotopy inverses.
It is not hard to see that
if $(\lqp,\rqp)$ is a Quillen equivalence, then both vertical arrows
are weak homotopy equivalences, and that these arrows are homotopy inverse to
one another in an appropriate sense \cite[17.5]{rDHKS}. In particular, if $\rqp$
preserves all equivalences, then the map
$\Moduli(\catMCp)\to\Moduli(\catMC)$ induced by $\rqp$ is a weak
homotopy equivalence.
\end{titled}

\begin{titled}{Derived mapping spaces}
If $X$ and
  $Y$ are objects of \catMC, $\maph_{\catMC}(X,Y)=\maph(X,Y)$ denotes the simplicial set
  of derived maps $X\to Y$; technically, this is defined in terms of
  the hammock localization \cite{calculating}, and it depends only on
  the equivalences in \catMC. If \catMC/ is a simplicial model
  category \cite[II.3]{goerssJardine}, then $\maph(X,Y)$ is canonically weakly homotopy equivalent to
  the  space $\Map_{\catMC}(X\cof,Y\fib)$  (see
  \cite{functionComplexes}).  The
  bifunctor $\maph({-},{-})$ converts an equivalence in either
  variable into a weak homotopy equivalence. By
  \cite[1.1]{calculating} the space $\maph(X,Y)$ is
  canonically weakly equivalent to the nerve of the category depicted
  pictorially as follows
   \begin{equation}\label{DiagramMapFormula}
   \xymatrix@R-14pt{     &  X_1 \ar@{.>}[dd]\ar[dl]_\sim\ar[r] 
                     &  Y_1\ar@{.>}[dd] \\
               X & {} & {} & Y \ar[ul]_\sim \ar[dl]^\sim \\
  &  X_2\ar[ul]^\sim \ar[r] &  Y_2}\raisebox{-27pt}{.}
  \end{equation}
  In this convention, zigzagging across the top
  gives one object of this category, zigzagging across the bottom
  another, and the entire commutative diagram with the dotted arrows
  drawn in represents a morphism from the top object to the bottom
  one. Dugger has shown \cite{duggerClassification} that if $X$ is cofibrant, $\maph(X,Y)$ is equivalent to the
  nerve of the less complicated category
\begin{equation}\label{SmallDiagramMapFormula}
\xymatrix@R-14pt { & Y_1\ar@{.>}[dd] & \\
         X\ar[ur]\ar[dr] && Y\ar[ul]_{\sim}\ar[dl]^{\sim}\\
            & Y_2 &}\raisebox{-28pt}{.}
\quad\quad
\end{equation}
\end{titled}

\begin{numbered}{Left proper model
    categories}\label{DiscussLeftProper} We need a slight
  generalization of Dugger's result. The model category 
 $\catMC$ is \emph{left proper} if the pushout of an equivalence along a
 cofibration is again an equivalence \cite[II.8.P2]{goerssJardine}. 
In \cite[2.7]{rezkProperModel}, Rezk observes that \catMC/ is left
 proper if and only if for any equivalence $Z\to Z'$ in \catMC/ the
 restriction functor $\rho:\Undercat {Z'}{\catMC}\to\Undercat Z{\catMC}$ is
the right adjoint of a Quillen equivalence.
Since $\rho$ preserves all weak equivalences, it follows that $\rho$ induces a weak equivalence
$\Moduli(\Undercat {Z'}{\catMC})\eqto\Moduli(\Undercat Z{\catMC})$
(see \ref{QPModuliSpaceDiagram} ff.).
  \begin{prop}
    Suppose that \catMC/ is left proper, and that $X$ and $Y$ are
    objects of \catMC/ such that $X\cof\mycoprod Y\to X\mycoprod Y$ is
    an equivalence. Then the nerve of the category
    \ref{SmallDiagramMapFormula} has the weak homotopy type of
    $\maph(X,Y)$.\end{prop}
\begin{proof}
Let $Z=X\cof\mycoprod Y$ and $Z'=X\mycoprod Y$.
 The nerve of
\ref{SmallDiagramMapFormula} is a union of components of
$\Moduli(\Undercat{Z'}{\catMC})$, and  according to Dugger, the
corresponding union of components of $\Moduli(\Undercat{Z}{\catMC})$
computes $\maph(X,Y)$. If $f:Z\to Z'$ is an
equivalence, it follows from Rezk's observation above that the map of
 moduli spaces induced by restriction over $f$ is a weak
equivalence. Since this map preserves the appropriate
components, the conclusion follows.
\end{proof}
\end{numbered}

\begin{rem}\label{DiscussModuliAndAut}
 Let $\Auth_{\catMC}(X)=\Auth(X)$ denote the
  union of those components of $\maph(X,X)$ which are
  invertible up to homotopy. There is a natural weak
    homotopy equivalence
    \begin{equation}\label{ModuliAndMaps}
        B\auth(X) \sim \Moduli(\catMC)_X\,,
    \end{equation}
 where $\Moduli(\catMC)_X$ is the component of the moduli space
 $\Moduli(\catMC)$ (\ref{DefineModuli})
 corresponding to~$X$. This is proved by stringing
  together results from \cite[5.5]{simplicialLocalizations},
  \cite[2.2]{calculating} and \cite[4.6(ii)]{functionComplexes}.
\end{rem}
\end{NumberedSubSection}


\begin{NumberedSubSection}{Quillen pairs and homotopy fibres of moduli
    spaces}\label{QuillenPairs}
 Suppose that $(\lqp,\rqp)$ is a Quillen pair as above, such that $\rqp$
preserves all equivalences and thus induces a map $\Moduli(\rqp)\colon\Moduli(\catMCp)\to
\Moduli(\catMC)$. We are interested in showing that the homotopy fibre of
$\Moduli(\rqp)$ over $A\in \Moduli(\catMC)$ is often given by
the nerve of the under category $\Undercat{A}{\Moduli(\rqp)}$.

\begin{prop}\label{ComputeHomotopyFibre}
   Suppose that $\lqp\colon\catMC\leftrightarrow\catMCp\colon\rqp$ is
   a Quillen pair such that $\rqp$ preserves all equivalences. Then
   the homotopy fibre of $\Moduli(\rqp)$ over~$A\in\Moduli(\catMC)$ is
   given by $\Undercat{A}{\Moduli(\rqp)}$ if either
   \begin{enumerate}
     \item $A$ is cofibrant, or
      \item \catMCp/ is left proper and $\lqp(A\cof)\to\lqp(A)$ is an
   equivalence.
   \end{enumerate}
\end{prop}

\begin{proof}
Consider the following two categories:
\[
\raisebox{-27.5pt}{$\Undercat{A}{\Moduli(\rqp)}:$}\hskip4pt\relax
\xymatrix@R-14pt{
     & \rqp X_1 \ar@{.>}[dd]_{\rqp f}\\
A\ar[ur]^\sim\ar[dr]_\sim\\
     & \rqp X_2
}
\xymatrix@R-14pt{{\vphantom{\rqp
      X_1}X_1}\ar@{.>}[dd]_\sim^f\\{\vphantom{A}}\\{\vphantom{\rqp X_2}X_2}
}
\qquad\
\raisebox{-27.5pt}{$\Moduli(\Undercat{\lqp A}{\catMCp}):$}\hskip4pt\relax
\xymatrix@R-14pt{
   & {\vphantom{\rqp X_1}X_1}\ar@{.>}[dd]_\sim\\
\lqp(A)\ar[ur]^{\alpha_1}\ar[dr]_{\alpha_2}\\
   &{\vphantom{\rqp X_2}X_2}
}
\]
The fact that $\rqp$ preserves equivalences
implies that $\Undercat{A}{\Moduli(\rqp)}$ is isomorphic to the  union of
the  components of $\Moduli(\Undercat{\lqp A}{\catMCp})$ containing
maps $\alpha:\lqp A\to X$ whose adjoint $\alpha^\flat:A\to\rqp X$ is
an equivalence. Two conclusions follow. First, if assumption (2) holds the
natural map $\Undercat{A}{\Moduli(\rqp)}\to
\Undercat{A\cof}{\Moduli(\rqp)}$ is a weak homotopy equivalence (\ref{DiscussLeftProper}), and hence without
loss of generality we can assume that (1) holds and $A$ is cofibrant. Secondly, if
$A\to A'$ is an equivalence between cofibrant objects of $\catMC$, the
natural map $\Undercat{A'}{\Moduli(\rqp)}\to
\Undercat{A}{\Moduli(\rqp)}$ is a weak homotopy equivalence \cite[2.5]{rezkProperModel}. 

  We now perturb the problem into
  one which can be solved by combining this last observation with
  Quillen's Theorem~B (\ref{QuillenTheoremB}). Let
  $\catMC\cof\subset\catMC$ be the subcategory of cofibrant objects,
  and consider the two categories described
  pictorially below (all arrows are equivalences).
\newcommand{\catLeft}{\Moduli(\catMCp)}
\newcommand{\catMiddle}{\Moduli_{\rqp}(\catMCp)}
\newcommand{\catRight}{\Undercat{\Moduli{(\catMC\cof)}}{\Moduli(\rqp)}}
\newcommand{\botLeft}{\Moduli(\catMC)}
\newcommand{\botMiddle}{\Moduli(\catMC)}
\newcommand{\botRight}{\Moduli(\catMC\cof)}
\[ \raisebox{-22pt}{$\catLeft\,:$}\quad
\xymatrix {{\vphantom{\rho(X_1)}X_1}\ar@{.>}[d]\\
  {\vphantom{\rho(X_2)}X_2}}
\qquad\qquad
\raisebox{-22pt}{$\catRight\,:$}\quad
  \xymatrix {{B_1\cof}\ar[r] \ar@{.>}[d] & \rqp(X_1)\ar@{.>}[d]_{\rqp(f)}
  \\
             {B_2\cof} \ar[r] & \rqp(X_2)} 
 \xymatrix{{\vphantom{\rho(X_1)}X_1}\ar@{.>}[d]_f \\{\vphantom{\rho(X_2)}X_2}}
           \]
The category $\catRight$ is a path space construction in which an 
object consists of
object $X_1\in\catMCp$, an object $B_1\cof\in\catMC\cof$, and an equivalence
$B_1\cof\to\rqp(X_1)$; a morphism consists of equivalences $X_1\to
X_2$ and $B_1\cof\to B_2\cof$ making the indicated diagram commute.
Let $u:\catMC\cof\to\catMC$ denote the inclusion. There is a diagram of functors
\[
    \xymatrix{{\catLeft}\ar[d]_{\Moduli(\rho)}  &
      {\catRight}\ar[l]_-v\ar[d]_w\\
    {\botLeft} &{\botRight}\ar[l]_{\Moduli(u)}}
\]
where $v$ picks out $X_1$ and $w$ picks out $B_1\cof$; the square
commutes up to an explicit natural transformation.
The functors $\Moduli(u)$ and $v$ are weak homotopy
equivalences, since there are oppositely oriented
functors (given by choosing functorial cofibrant replacements) such that composites are connected to appropriate identity functors
by natural transformations. 
If $A$ is cofibrant, it is easy to construct
a functor $\Undercat{A}{\Moduli(\rqp)}\to\Undercat{A}{w}$
(natural in~$A$) which is a weak homotopy equivalence. As above, this
implies that any equivalence $A\to A'$ in $\catMC\cof$ induces a weak
homotopy equivalence $\Undercat{A'}{w}\to\Undercat{A}{w}$.  By Quillen's
Theorem~B, the homotopy fibre of $w$ over $A$, or equivalently the
homotopy fibre of $\Moduli(\rqp)$ over~$A$, is weakly homotopy
equivalent to $\Undercat Aw\sim\Undercat A{\Moduli(\rqp})$.
\end{proof}
\end{NumberedSubSection}

\begin{NumberedSubSection}{Quillen pairs and mapping spaces}
We will also need that Quillen pairs are topologically adjoint with
respect to model-category theoretic mapping spaces.
A version of this up to homotopy is included in
\cite[5.6.2]{hovey} (cf \cite[17.4.15]{hirschhorn}), but we prefer a slightly more rigid formulation.
As usual, $A\cof\eqto A$ and $Y\eqto Y\fib$ are cofibrant and
fibrant replacements.

\begin{thm}\label{AdjointFunctionComplexes}
  Suppose that $\lqp\colon\catMC \leftrightarrow\catMCp\colon\rqp$ is
  a Quillen pair, and that $A$ and $Y$ are objects of $\catMC$ and
  \catMCp/ respectively. Then there is a natural weak homotopy
  equivalence
  \[
       \maph_{\catMC}(A,\,\rqp(Y\fib)) \sim \maph_{\catMCp}(\lqp(A\cof),\,
       Y)\,.
  \]
\end{thm}

\begin{proof}
We can assume without loss of generality
that $A$ is cofibrant and $Y$ is fibrant. 
Consider the category $\SquigThree$ described by the following diagram:
\begin{equation}\label{SquigThreeDisplay}
\xymatrix@R-14pt{     &  A_1 \ar@{.>}[dd]\ar[dl]_\sim\ar@{~>}[r] 
                     &  Y_1\ar@{.>}[dd] \\
               A & {} & {} & Y \ar[ul]_\sim \ar[dl]^\sim \\
  &  A_2\ar[ul]^\sim \ar@{~>}[r] &  Y_2}
\end{equation}
where an undulating arrow $\xymatrix{A_i\ar@{~>}[r]& Y_i}$ represents
(equivalently) a map
$A_i\to\rqp(Y_i)$ in \catMC/ or a map $\lqp(A_i)\to Y_i$ in \catMCp.
We will show that $\SquigThree$  is weakly equivalent to $\maph(\lqp A, Y)$; a
dual argument shows that it is weakly equivalent to $\maph(A,\rqp Y)$.

Let $\catMC\cof$ be the category of cofibrant objects in~\catMC.
A functorial factorization argument shows that $\SquigThree$ is weakly
homotopy equivalent to the full subcategory $\SquigThree\cof$
consisting of zigzags with $A_1\in\catMC\cof$. Picking off
$A_1\eqto A$ gives a functor
$F\colon\SquigThree\cof\to\Overcat{A}{\Moduli(\catMC\cof)}$.
Given an object $U=\langle A'\eqto A\rangle$ of $\Overcat{A}{\Moduli(\catMC\cof)}$,
the over category $\Undercat{U}{F}$ has objects consisting of diagrams
of the form
\[\xymatrix@R-18pt{
   & A'\ar[dl]_\sim\ar[dd]      \\
A                        \\
   &A_1\ar[ul]^\sim \ar@{~>}[r] &Y_1 & Y \ar[l]^\sim
}
\]
This is an \emph{object} of $\Undercat{U}{F}$; in this category, $A$,
$A'$ and $Y$ are fixed, but $A_1$ and $Y_1$ are allowed to vary. This
 category is  homotopy equivalent to the subcategory
consisting of objects in which $A'\to A_1$ is the identity map; by
\ref{SmallDiagramMapFormula}, this subcategory has the weak homotopy type of
$\maph_{\catMCp}(\lqp (A'), Y)$ and so in particular the weak homotopy
type is independent of the choice $U$ of object. By Quillen's Theorem~B,
$\maph_{\catMCp}(\lqp A,Y)$ is equivalent to the homotopy fibre of
$F$. Since the target $\Overcat{A}{\Moduli(\catMC\cof)}$ of~$F$ is contractible (this
category has a terminal object), it follows that
$\maph_{\catMCp}(\lqp A,Y)$ is weakly equivalent to $\SquigThree$.
\end{proof}

\end{NumberedSubSection}

\section{Monoidal categories and the proof of \ref{GeneralTheorem}}\label{MonoidalCategories} In this section we work out the 
machinery sketched in \ref{Connection}. The triple
$(\catM,\gentensor,\genunit)$ will be a monoidal category in the sense
of \cite[VII.1]{maclane}, and \catMon/ the category of monoids
 in \catM/ \cite[VII.3]{maclane}. For the rest of this section
$R$ and $S$ are two fixed objects of \catMon. We use the following
notation:
\begin{itemize}
\item \catMods/ $=$ the category of right $S$-modules,
\item \catModrs/ $=$ the category of \bimod RSs, and
\item \catMonrs/ $=$ the category of monoidal objects under $R$ and $S$.
\end{itemize}
A \emph{point} for an object $X\in
\catM$ is defined to be a \catM-map $\genunit\to X$; a map between pointed
objects is required to respect the points. A superscript $\myunit$
indicates a category of pointed objects, so that for instance
\begin{itemize}
\item \catBirsp/ $=$ the category of pointed \bimod RSs.
\end{itemize}
Note that the category of pointed \bimod RSs is isomorphic to the
category of \bimod RSs under $R\gentensor S$. We now introduce the
axioms under which \ref{GeneralTheorem} holds.

\begin{axiomone}\label{AxiomOne} The categories \catM, \catMon,  \catMods,and
  \catModrs/ possess compatible (see below) model category structures.
\end{axiomone}

The term \emph{compatible} means that a map in one of these
categories is an equivalence or a fibration if and only if the
underlying map in \catM/ has the same property. For general reasons,
the above model category structures on \catMon/ and on \catBirs/ 
extend  to compatible
model structures on \catMonrs/ and \catBirsp/ \cite[1.1.8$\,${ff}]{hovey}.

There is a forgetful functor 
$\forget\colon \catMonrs\to\catBirsp$, with the point
for $\forget(T)$  provided by the unit map $\genunit \to T$.

\begin{axiomtwo}\label{AxiomTwo} The functor $\forget\colon \catMonrs\to\catBirsp$
  has a left adjoint $\Envel$.
\end{axiomtwo}

\begin{rem}\label{AutoQuillenPair}
  The compatibility condition in \ref{AxiomOne}
  guarantees that $\forget$
  preserves fibrations and 
  equivalences, and it follows  that $(\Envel,\forget)$ forms a Quillen pair.
  We will let $\Envelh$ denote the left derived functor of $\Envel$,
  which  is
  given  by $\Envelh(X)=\Envel(X\cof)$, where $X\cof\to X$
  is a functorial cofibrant replacement in \catBirsp. See
  \cite[1.3.2]{hovey}, but note that for us the codomain of
  $\Envelh$ is \catMonrs, not the homotopy category of~\catMonrs.
\end{rem}
Given $X\in\catModsp$,
there is a right $S$-module map
\[
    S\cong \genunit\gentensor S\to X
\]
derived from the point in $X$. The object $X$ is
\emph{distinguished} if this map is an equivalence. An object of
\catBirsp/ is distinguished if it is distinguished as a pointed right
$S$-module. Similarly, an object  $T\in\catMonrs$ is distinguished if the
structure map $S\to T$ is an equivalence.

\begin{axiomthree}\label{AxiomThree} The derived functor
  $\Envelh:\catBirsp\to\catMonrs$ preserves distinguished objects.
\end{axiomthree}

The remaining axioms are more technical and have the flavor of
non-degeneracy assumptions.
Let $R\cof\to R$ be a cofibrant replacement for $R$ in \catMon.

\begin{conditiond}\label{AxiomFour}
 One of the following two conditions holds:
 \begin{enumerate}
 \item $R$ itself is cofibrant as an object of \catMon, or
 \item \catMon/ is left proper (\ref{DiscussLeftProper}), and the map
   $R\cof\mycoprod^{\catMon} S\to R\mycoprod^{\catMon} S$ is an equivalence.
 \end{enumerate}
\end{conditiond}

In the next statement, $S\cof\to S$ is a cofibrant replacement in~\catMods.

\begin{conditionf}\label{AxiomFive}
One of the following two conditions holds:
\begin{enumerate}
\item $S$ itself is cofibrant as an object of \catMods, or
\item \catModrs/ is left proper, and the map
  $R\gentensor S\cof\to R\gentensor S$ is an equivalence.
\end{enumerate}
\end{conditionf}

\begin{rem}\label{GeneratorCofibrant}
  It is very tempting to assume that $S$ is necessarily 
 cofibrant as an object of \catMods,
  since $S$ is the free right $S$-module on one generator. But
  the notion of ``one generator'' here is tricky: it is more accurate
  to say that $S$ is the
  free right $S$ module on the object \genunit/ of \catM. The real
  issue is whether \genunit/ is
  cofibrant in \catM.
\end{rem}

 If $M$ is an \bimod RR and $\monmap:R\to S$ is a map of monoidal
objects,  write $M\gentensor^{\monmap}_RS$ for the
coequalizer of the two maps $M\gentensor R\gentensor S\to M\gentensor
S$ obtained by pairing the central $R$ either with $M$  or
(via~$\monmap$) with~$S$; the coequalizer is to be computed in the category
\catBirs. The functor $M\mapsto M\gentensor_R S$ is left adjoint
to the functor $\catBirs\to\catBirr$ induced by composition with~$\monmap$.

\begin{axiomsix}\label{AxiomSix} Suppose that $\monmap:R\to S$ is a map of
  monoidal objects, $\genunit\cof\to\genunit$ is a cofibrant
  replacement for $\genunit$ in \catM, and $R\cof\to R$ is a cofibrant
  replacement for $R$ in \catBirr. Then the following two conditions
  hold.
  \begin{enumerate}
  \item $\genunit\cof\gentensor S\to \genunit\gentensor S\cong S$ is
    an equivalence (in \catMods), and
  \item $R\cof\gentensor^{\monmap}_{R}S\to R\gentensor^{\monmap}_RS\cong S$ is an
    equivalence (in \catBirs). 
  \end{enumerate}
\end{axiomsix}

We first prove \ref{DeloopGeneralTheorem}, which relates maps between
monoid objects to moduli spaces of pointed bimodules.
Recall from \ref{ModelGeneralities} that  $\maph_{\catMon}(R,S)$
denotes the space of maps $R\to S$ provided by the model category
structure on \catMon. 
Let \catBirspd/ be the subcategory of \catBirsp/ given by
the distinguished objects; all morphisms in this category are
equivalences, so
$\Moduli(\catBirspd)=\catBirspd$ is just the nerve.

\begin{thm}\label{IdentifyMonoidMaps}
  If axioms I-IV hold, then 
  there is a natural weak equivalence of spaces
  \[         \maph_{\catMon}(R,S)\sim \Moduli(\catBirspd)\,.\]
\end{thm}

\begin{titled}{Proof of \ref{IdentifyMonoidMaps}}
Combining \ref{AxiomFour} with \ref{DiscussLeftProper} shows
 there is a natural weak
equivalence $\maph_{\catMon}(R,S)\sim\Moduli(\catMonrsd)$, where 
\catMonrsd/ is the category of all distinguished objects in \catMonrs/
and all maps (necessarily equivalences) between them. Recall that for
$X\in\catBirsp$, $\Envelh(X)$ can be computed as $\Envel(X\cof)$, where
$X\cof\RightArrow{\sim}X$ is a functorial cofibrant replacement
for~$X$. The functor $\forget$ restricts to a functor
$\catMonrsd\to\catBirspd$ and it follows from \ref{AxiomThree} that
$\Envelh$ restricts to a functor $\catBirspd\to\catMonrsd$. The arrows
\[
   \Envel((\forget T)\cof)\to T
\quad\text {and}\quad
   X \leftarrow X\cof \to \forget (\Envel(X\cof))
\]
show that the composites $\Envel\forget$ and $\forget\Envel$ of these
restricted functors  are each
connected to the respective identity functor by a chain of natural
transformations; it follows \catBirspd/ and \catMonrsd/ have weakly
equivalent nerves. \qed
\end{titled}

Call a right $S$ module \emph{potentially distinguished} if it is
(abstractly) equivalent to $S$ itself, and let \catModsd/ denote the full
subcategory of \catMods/ containing the objects which are potentially
distinguished.  Similarly, let \catBirsd/
be the full subcategory of \catBirs/ containing bimodules which are
potentially distinguished as right $S$-modules.

The next statement unbundles the basepoint from \ref{IdentifyMonoidMaps}.

\begin{thm}\label{MonoidMapFibration}
  If axioms I-V hold, then there is a natural fibration
  sequence
  \[  \maph_{\catMon}(R,S)\to \Moduli(\catBirsd)\to\Moduli(\catModsd)\,.
  \]
  where the fibre is to be taken over $S\in\catModsd$.
\end{thm}

\begin{proof}
  Consider the adjoint functors
  $\lqp\colon\catMods\leftrightarrow\catBirs\colon\rqp$, where $\rqp$
  forgets the left $R$-structure and $\lqp(X)=R\gentensor X$. As in
  \ref{AutoQuillenPair},
  $(\lqp,\rqp)$ forms a Quillen pair. By \ref{AxiomFive} and
  \ref{ComputeHomotopyFibre}, the homotopy fibre of $\Moduli(\rqp)$
  over $S$ is naturally weakly homotopy equivalent to
  $\Moduli(\catBirspd)$, and hence, by \ref{IdentifyMonoidMaps} to
  $\maph_{\catMon}(R,S)$. The theorem follows from the fact that
  $\Moduli(\catModsd)$ is the component of $\Moduli(\catMods)$
  containing $S$, while $\Moduli(\catBirsd)$ is the inverse image of
  this component in $\Moduli(\catBirs)$.
\end{proof}

Finally, we loop down the fibration sequence from
\ref{MonoidMapFibration} and rewrite the spaces involved. 

\begin{thm}\label{LoopedMonoidMapFibration}
Suppose that axioms I-VI hold, and that $\monmap:R\to S$ is a map of
monoidal objects. Let $S_{\monmap}$ denote $S$ considered via~$\monmap$ as an \bimod
RR. Then there is a natural fibration sequence
\[
  \Omega\maph_{\catMon}(R,S)_{\monmap}\to 
        \maph_{\catBirr}(R,S_{\monmap})\to\maph_{\catM}(\genunit,S)\,.
\]
\end{thm}

\begin{rem}
  In this statement,  $\Omega\maph_{\catMon}(R,S)_{\monmap}$ denotes the loop
  space taken with~$\monmap$ as the basepoint, and $R$ is to be 
  treated as an
  \bimod RR in the natural way. The right hand map is induced by the
  unit $\genunit\to R$, and the fibre is meant to be taken over the
  unit $\genunit\to S$.  
\end{rem}

\begin{titled}{Proof of \ref{LoopedMonoidMapFibration}}
 Looping down the fibration sequence from \ref{MonoidMapFibration}
 gives a sequence
\[
   \Omega\maph_{\catMon}(R,S)_{\monmap} \to \Omega\Moduli(\catBirsd)_{\monmap}
   \to\Omega\Moduli(\catModsd)_S\,.
\]
We begin by considering the right hand space. By general properties of
moduli spaces (\ref{ModuliAndMaps}), this loop space is equivalent to the
subspace $\auth_{\catMods}(S)$ of $\maph_{\catMods}(S,S)$ consisting
of homotopically  invertible maps.  
By \ref{AxiomSix}(1) and \ref {AdjointFunctionComplexes}
(this last applied to the forgetful functor $\catMods\to \catM$ and
its left adjoint ${-}\gentensor S$) the space $\maph_{\catMods}(S,S)$
is weakly equivalent to $\maph_{\catM}(\genunit, S)$. 

For similar reasons, the middle space  is equivalent to the
subspace $\auth_{\catBirs}(S_{\monmap})$ of $\maph_{\catBirs}(S_{\monmap}, S_{\monmap})$ consisting of homotopically
invertible maps. By \ref{AxiomSix}(2) and
\ref{AdjointFunctionComplexes} (this last applied to the restriction
functor $\monmap^*:\catBirs\to\catBirr$ and its left adjoint
${-}\gentensor^{\monmap}_RS$), the space $\maph_{\catBirs}(S_{\monmap}, S_{\monmap})$ is naturally weakly
equivalent to $\maph_{\catBirr}(R,S_{\monmap})$. 

All in all, there is a commutative diagram
\[
\xymatrix{
 {\auth_{\catBirs}(S_{\monmap})} \ar[d]\ar[r]
                                &
 {\maph_{\catBirs}(S_{\monmap},S_{\monmap})}\ar[d]\ar@{<->}[r]^\sim &{\maph_{\catBirr}(R,S_{\monmap})}\ar[d]\\
  {\auth_{\catMods}(S) \ar[r]} &
  {\maph_{\catMods}(S,S)}\ar@{<->}[r]^\sim &{\maph_{\catM}(\genunit, S)}
 }
\]
in which the left horizontal arrows are component inclusions. Since a map
$S_{\monmap}\to S_{\monmap}$ in \catBirs/ is an equivalence if and only if the
underlying map of right $S$-modules is an equivalence, the left hand
square  is 
a homotopy fibre square. Consequently, the homotopy fibre over the
identity map of $S$ on the left, namely $\Omega\maph_{\catMon}(R,S)_{\monmap}$,
is equivalent to the homotopy fibre over the image of this identity
map on the right.  We leave it to the reader to check that this image
is the unit $\genunit\to S$, and that the right vertical map is as
described. \qed
\end{titled}

\section{Model structures, free operads and coproducts of operads}
  \label{CModelFreeCoprod}

\renewcommand{\catM}{\nspace}
\renewcommand{\catMon}{\operads}
In this section we  develop a few homotopical properties of operads. Deep
in the background is the symmetric monoidal category
$(\sspace,\times,*)$ of simplicial sets with cartesian product, but
the ambient monoidal category in this section is the category $\nspaces=
(\nspaces,\ctensor,\cunit)$ of \gradedspaces, with the (nonsymmetric)
composition monoidal structure from \ref{MoreDetails}.  An
\emph{operad} is a monoid in $\nspaces$.
The category \sspaces/ has a simplicial model category structure in
which the equivalences are the weak homotopy
equivalences and the cofibrations are the monomorphisms
\cite[I.11]{goerssJardine}.
This model structure extends to
\nspaces/ by declaring a map in \nspaces/ to be an
equivalence (resp. cofibration, fibration) if and only if on each
level  it gives an equivalence (resp. cofibration, fibration) in
\sspaces.

Let 
$\catMp$ denote the category of \gradedspaces/ furnished with a
basepoint at level~one, let \catMon/ denote the category of operads, and for two chosen operads
$\opa,\opb$, let
\catModb/ and \catModab/ denote respectively the
categories of  right \opb-modules and \bimod\opa\opb s.   The
forgetful functor $\catMon\to\catM$ has a left adjoint $\freeop$ (the
free operad functor), and the
forgetful functor $\catMon\to\catMp$ obtained by retaining the
unit as a basepoint has a left adjoint $\freeopp$. In this section we
study model structures on these categories, and prove that  $\freeop$,
$\freeopp$, and the coproduct functor on $\catMon$ are homotopy invariant.

\begin{prop}\label{RezkModelTheorem}\cite[\S3]{rezkThesis} \cite{rezkProperModel} The categories \catMon, 
  \catModb, and \catModab/  have simplicial model structures compatible
  (\ref{AxiomOne}) with
  the model structure on \nspace.
\end{prop}

\begin{rem}\label{DescribeCofibrations}
 Rezk produces these model category
  structures in the setting of $\Sigma$-operads, but the case of
  planar operads is a bit simpler; the arguments follow the lines of
  \cite[II.4.1]{goerssJardine} and \cite[II.5.1]{goerssJardine} (see
  also \cite[7.1]{rezkProperModel}).   It will be useful later on to
  understand the cofibrations in these model categories. As in
  \cite[\S6]{rezkProperModel}, a \emph{degeneracy object} is a
  simplicial object without face operators.
   According to
  \cite[\S6]{rezkProperModel}, a map $X\to Y$ in one of these model
  categories is a cofibration if it is (a retract of) a monomorphism with
  the  additional property that, as a degeneracy object, $Y$
  is isomorphic to the coproduct of $X$ with a free degeneracy object.
  The notions of \emph{free} and \emph{coproduct} here are to be interpreted
  in the relevant operad or (bi)module setting.
\end{rem}

\begin{prop}\label{CoproductInvariance}
  The coproduct construction on operads preserves equivalences. The
  category \catMon/ is a left proper model  category.
\end{prop}

\begin{prop}\label{PointedFreeInvariance}
  The functors $\freeop$ and $\freeopp$ preserve equivalences.
\end{prop}

The rest of this section is devoted to proofs of
\ref{CoproductInvariance} and \ref{PointedFreeInvariance}; these are routine but
involve a substantial amount of notation and bookkeeping.

\begin{numbered}{Identities and compositions}
  The identity element of an operad \opa/ is denoted~$\Identity\opa$.
  Given a point $a\in\opa\Of n$, and $n$ elements
  $b_i\in\opa\Of{k_i}$ ($1\le i\le n$), we use composition
  notation $a(b_1,\ldots,b_n)$
  to denote the image of a
  tuple $(a,b_1,\ldots,b_n)$ under the operad structure map.
\end{numbered}

\begin{numbered}{Trees and free operads}
We first describe the set of isomorphism classes of (planar rooted)
trees. Each tree
$\tree$ has a height $\height(\tree)$, a leaf count $\lc(\tree)$, and a
\branchcount/~$\bc(\tree)$.
The set of trees is defined
recursively on height by declaring that there is  one tree of height~0, called the
  trivial tree and denoted \etree. Its leaf
  count is~1, and its \branchcount/ is~0.
A tree of height $k>0$ is then a tuple
  $\treeof(n;\tree_1,\ldots,\tree_n)$, where $n>0$ is an integer, each
  $\tree_i$ is a tree of height less than~$k$, and at least one
  $\tree_i$ has height~$k-1$. 
The numerical invariants of a nontrivial tree
are given inductively
  by 
 \[
  \begin{aligned}
    \height\treeof(n;\tree_1,\ldots,\tree_n) &= 1 +
    \operatorname{max}\{\height(\tree_i)\} \text{\quad\qquad height}\\
    \lc\treeof(n;\tree_1,\ldots,\tree_n) &= \sum \lc(\tree_i)
    \text{\quad\quad\qquad\qquad leaf count}\\
    \bc\treeof(n;\tree_1,\ldots,\tree_n) &= 1 + \sum\bc(\tree_i)
    \text{\quad\ \quad\qquad \branchcount/.}
  \end{aligned}
  \]
Here are some examples.
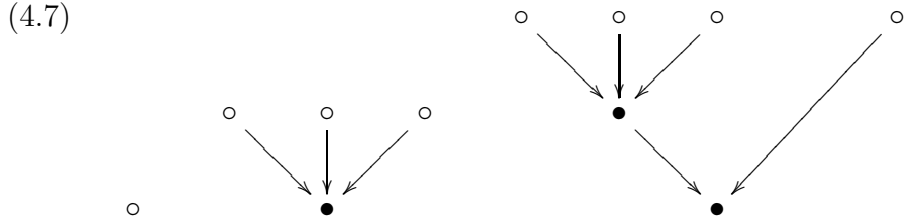
\begin{equation}\label{PicturesOfTrees}
\xymatrix{{\vphantom{\circ}}\\ {\vphantom{\circ}}\\ \circ}
\qquad
\xymatrix{{\vphantom{\circ}}\\
\circ\ar[dr] & \circ \ar[d]       & \circ\ar[dl] \\
             {}       & \bullet  }
\qquad
\xymatrix{\circ\ar[dr] & \circ \ar[d]       & \circ\ar[dl]  & {} & \circ\ar[ddll]\\
              {}       & \bullet \ar[dr]    \\
              {}       & {}                 &\bullet
}
\end{equation}
On the left is \etree: one leaf resting on the ground, no
branching. The middle tree is $\treeof(3;\etree,\etree,\etree)$: height
one, three leaves, one branching. The one on the right is
$\treeof(2;{\treeof(3;\etree,\etree,\etree)},\etree)$: height 2,
4~leaves, 2~branchings. Trees of branching at most~$1$, such as the
trees on the left above, are sometimes called \emph{corollas}.

Given \agradedspace/ $X$, define a space $\treelab X{\tree}$ for each
tree~$\tree$ by declaring $\treelab X{\etree} = *$, and, if
$\tree=\treeof(n;\tree_1,\ldots,\tree_n)$, setting
\[ \treelab
  X{\tree}=X\Of n\times{\textstyle \prod_i \treelab X{\tree_i}}\,.
\]
This is called the \emph{space of labelings of the tree~$\tree$
  with labels from~$X$}.  We will denote such a labeling~$l$ by
  $\treeof(x;l_1,\ldots,l_n)$, where $x\in X\Of n$ and $l_i$ is a
  labeling of $\tree_i$, i.e, $l_i\in\tree_i(X)$. We treat $\treelab X{\tree}$ as \agradedspace/
concentrated at level~$\lc(\tree)$, and  let $\Treelab X$
denote the union of the \gradedspaces/ $\treelab X\tree$, taken over all
trees~$t$. This is the space of \emph{trees with labels from~$X$}.
In the pictorial terms of \ref{PicturesOfTrees}, a labeling of a tree
by $X$ is a choice, for each solid dot in the sketch of the tree, of
a simplex in $X_n$, where $n$ is the number of edges pointing inward
towards the dot.

\begin{rem}\label{LabeledTreeNotation}
  For $x\in X\Of n$, it is convenient to use $\genof x$ to denote the
  labeling $\treeof(x;l_1,\ldots,l_n)$ in which each $l_i$ is the
  unique labeling of the trivial tree. The map $x\mapsto\genof x$
  gives a map $X\to\Tof (X)$ of \gradedspaces.  Each labeling $l$ in
  $\Tof(X)$ has an underlying tree $\shape l$, called the \emph{shape}
  of~$l$; the labeling inherits its depth $\height(l)$, leaf count
  $\lc(l)$ and \branchcount/ $\bc(l)$ from the shape. The assignment
  $l\mapsto\shape l$ gives a map from $\Tof(X)$ to the graded discrete
  space $\Trees$ of all trees.  The \emph{set $\labset l$ of labels
    in~$l$} is defined by initializing $\labset *=\emptyset$ and
  constructing $\labset{\treeof(x;l_1,\ldots,l_n)}$ by adjoining $x$
  to the union of the sets $\labset {l_i}$, $1\le i\le n$. This object
  is a set of simplices in~$X$. The assignment $l\mapsto\labset l$
  gives a map from $\Tof (X)$ to the powerset of the union $\cup_nX_n$.
\end{rem}

The following proposition is elementary.

\begin{prop}\label{FreeOperadDescription}
  If $X$ is \agradedspace,  $\Treelab X$ is
  naturally isomorphic to the free operad $\freeop(X)$.
\end{prop}

\begin{rem}
  The operad composition in $\Tof(X) $ is defined as follows.
  Given a labeling $a\in \Tof 
  (X)$ with $\lc(a)=m$ and $m$ labelings $b_i\in \Tof(X)$, $1\le i\le m$, it is
  sufficient to specify the composite labeling $a(b_1,\ldots,b_m)$. 
This  is given
  inductively by setting  $*(b_1) = b_1$ and specifying that
  if $a=\treeof(x;l_1,\ldots,l_n)$ then
  \[
    a(b_1,\ldots,b_m)=
    \treeof(x;
       {l_1(b_1,\ldots,b_{\lc(l_1)}),\ldots,
       l_n(b_{m-\lc(l_n)+1},\ldots,b_m)})\,.
  \]
  The unique labeling of the trivial tree~\etree/ serves as
  $\Identity{\Tof(X)}$.
\end{rem}

\begin{numbered}{Reflexive coequalizers}\label{ReflexiveCoequalizers}
A \emph{reflexive pair} in some category
is a diagram which looks like the
one-truncation of a simplicial object; 
more explicitly, it consists of two objects $X_1$, $X_0$, maps
$d_0,d_1\colon X_1\to X_0$, and a map $s_0:X_0\to X_1$ such that both
composites $d_0s_0$ and $d_1s_0$ give the identity map of~$X_0$.
An augmentation for the pair is an object $X$ and a map $d:X_0\to X$ such that $dd_0=dd_1$.
The colimit of a reflexive pair is isomorphic to the coequalizer of
$(d_0,d_1)$. An augmented reflexive pair is \emph{exact} if the
natural map from the coequalizer of $(d_0,d_1)$ to $X$ is an
isomorphism, in which case the diagram is called a reflexive
coequalizer diagram. 
Reflexive coequalizers commute with finite products in the category of
sets \cite[3.2]{rezkProperModel}, and it follows
easily from this that \emph{reflexive coequalizers of just about any
algebraic species can be computed in the appropriate underlying
category}. See  \cite[4.3]{rezkProperModel} or the monograph \cite{rARV}.
In particular, {reflexive coequalizers of operads of spaces can be computed in
the category of \gradedspaces}.
\end{numbered}

 \begin{numbered}{A formula for the coproduct}
 Let $U:\catMon\to\nspaces$ be the forgetful functor, i.e.,
the right adjoint to~$\freeop$, and let $\comp\colon\catMon\to \catMon$
be the composite $\freeop U$. For any operad \opa, adjointness
considerations give
maps $\eval_{\opa}\colon \comp(\opa)\to \opa$ and
$\sigma_{\opa}:\comp(\opa)\to\comp^2(\opa)$, and it is not hard to check that
these give rise to a natural exact augmented reflexive pair
$\comp^2(\opa)\ReflexArrow\comp(\opa)\to\opa$ with $d=\eval_{\opa}$, $d_0=\eval_{\comp\opa}$,
$d_1=\comp(\eval_{\opa})$, and $s_0=\sigma_{\opa}$.
Let $\seqcoprod$ denote coproduct on the category of
\gradedspaces. Since colimits commute with one another and $\freeop$
as a left adjoint 
commutes with colimits, it follows that for any two operads \opa,
\opb/  there is a natural right exact
  augmented reflexive pair 
\begin{equation}\label{CoproductCoequalizer}
    \freeop(\comp\opa\seqcoprod\comp\opb)\ReflexArrow
    \freeop(\opa\seqcoprod \opb)\to \opa\mycoprod\opb\,.
\end{equation}
We have suppressed forgetful functors here (e.g., $\opa\seqcoprod\opb$
really signifies $U\opa\seqcoprod U\opb$) and in order to avoid
clutter we will continue to do this in what follows.
For the rest of this section we fix \opa/ and \opb, 
set $X_1=\freeop(\comp\opa\seqcoprod\comp\opb)$,
$X_0=\freeop(\opa\seqcoprod \opb)$, and let  $d_0$, $d_1$ and $s_0$ 
refer to the maps in the above reflexive pair involving $X_1$ and $X_0$.
Note again that the operad coproduct 
$\opa\mycoprod\opb$ is isomorphic as \agradedspace/ to the
coequalizer of  $(d_0,d_1)$ \emph{in the category of
  \gradedspaces/}, because this coequalizer is reflexive.

The simplices of $X_0$ are
labeled trees $l=\treeof(a;l_1,\ldots,l_n)$ in which the labels, e.g.,
$a$, range over simplices of \opa/ and of~\opb; simplices of
$X_1$ are labeled trees $v=\treeof(u;v_1,\ldots,v_n)$ in which
the labels are taken from $\comp(\opa)$ and $\comp(\opb)$.
The maps
$d_0$, $d_1$ and $s_0$ preserve identities; their inductive descriptions
are
  \begin{equation}\label{FormulaForCoproductDs}
   \begin{aligned}
    d_0\treeof(u;v_1,\ldots,v_n) &= u(d_0v_1,\ldots, d_0v_n)\\
    d_1\treeof(u;v_1,\ldots,v_n) &=
    \treeof({\eval(u)};{d_1(v_1)},\ldots,{d_1(v_n)})\\
    s_0\treeof(a;l_1,\ldots,l_n) &= \treeof(\genof
    a;{s_0(l_1),\ldots,s_0(l_n)})\,. 
    \end{aligned}
    \end{equation}
  In each case the label $u$ is a simplex of
  $\comp\opa\seqcoprod\comp\opb$  in level~$n$, treated in the obvious way as a
  simplex of $\freeop(\opa\seqcoprod\opb)$; $a$~is a simplex of
  $\opa\seqcoprod\opb$.  In the upper formula the
  composition operation takes place in $\freeop(\opa\seqcoprod\opb)$. In
  the middle formula the composition operations implicitly involved in
  computing $\eval(u)$ take place either 
  in $\opa$ or in $\opb$ (depending on
  whether $u\in\comp\opa$ or $u\in\comp\opb$). It is clear from the
  lower formula that a simplex~$v$  in~$X_1$ is in the image of~$s_0$ if
  and only if each label in~$v$ has \branchcount/~$1$.\end{numbered}

\end{numbered}

\begin{numbered}{An inductive construction of the coproduct}
  It is helpful to introduce  some terminology.

  \begin{defn}\label{DefinitionOfFour} 
  Suppose that \opaa/ and \opbb/ are operads, and 
 that $l=\treeof(a;l_1,\ldots,l_n)$  is a nonidentity simplex
    of $\freeop(\opaa\seqcoprod\opbb)$.
    \begin{enumerate}
      \item $l$ has \emph{color \opaa/ (resp. color \opbb)} if $a\in\opaa$
      (resp. $a\in\opbb$),
      \item $l$ \emph{displays adjacent colors} if any $l_i$~does, or
        if any  $l_i\neq\etree$ has the same color
      as $l$ ($\etree$ itself does not display adjacent colors),
     \item $l$ is \emph{unital} if $\Identity\opaa\in\labset l$ or
      $\Identity\opbb\in\labset l$, and
     \item  $l$ is \emph{collapsible}
      if $l$ displays adjacent colors or $l$ is unital.
   \end{enumerate}
  \end{defn}
Let $\Treelabpair\opa\opb k\subset \Tof(\opa\gradedspacesum\opb)=X_0$
denote the \gradedspace/ of all labeled trees of \branchcount/~$k$, and
$\collTreelab\opa\opb k$ the subobject of collapsible labeled trees.
Let $\Filt k$ be the subspace of $\opa\mycoprod\opb$ given by the union
of the images of $\Treelabpair\opa\opb j$ for $j\le k$.
The main proposition we use to prove homotopy invariance of
coproducts is the following one. 

  \begin{prop}\label{CoprodPushout}
    For each $k\ge 1$ there is a pushout diagram of \gradedspaces/
  \begin{equation*}
    \xymatrix{
           \collTreelab\opa\opb{k} \ar[r]\ar[d] &
                       \Treelabpair\opa\opb{k}\ar[d]\\
           \Filt{k-1}\ar[r] & \Filt{k}\
           }
    \end{equation*}
    which is also a homotopy pushout diagram.
  \end{prop}

\begin{titled}{Proof of \ref{CoproductInvariance} (given
    \ref{CoprodPushout})}
The \gradedspace/  $\Treelabpair\opa\opb k$ is homotopy invariant as a
functor of \opa/ and \opb, since it is formed by taking disjoint
unions of products of the graded constituents of these two
operads. The subobject $\collTreelab\opa\opb{k}$ is similarly homotopy
invariant, since it is determined by selecting appropriate components (picking labeled
trees that display adjacent colors) or by imposing basepoint
restrictions on appropriate cartesian factors (picking unital labeled
trees). Since $\Filt0=\{\Identity{\opa\mycoprod\opb}\}$ is a
single point, it follows by induction that each object $\Filt k$
depends in a homotopy invariant way on \opa/ and \opb, and so
$\opa\mycoprod\opb=\colim_k\Filt k$ also does. The statement that
\operads/ is left proper follows as in \cite[9.1]{rezkProperModel}
from the fact that the coproduct construction is homotopy invariant.
\qed
\end{titled}

In order to obtain \ref{CoprodPushout} from the coequalizer
presentation \ref{CoproductCoequalizer}, we will have to trim
$X_1$ a bit, without changing the
value of the \gradedspace/ coequalizer~\ref{CoproductCoequalizer}.
Note that the terminology of \ref{DefinitionOfFour} applies to
$X_1=\Tof(\comp \opa\gradedspacesum\comp\opb)$; moreover, the labels with
which the trees in $X_1$ are decorated lie in $\comp\opa$ or
$\comp\opb$, and thus have \branchcounts/ in their own right. In the
context of \ref{DefinitionOfFour}, say that a labeled tree is
\emph{alternating} if it does not display adjacent colors. 
Define two subobjects $\longX$ and
$\shortX$ of $X_1$ as follows:
\begin{itemize}
\item $v\in \longX$ if $v$ is  alternating, $\Identity{\comp\opa}$ and
  $\Identity{\comp\opb}$ do not appear as labels in $v$, and at least
  one label in~$v$  has
  \branchcount~$>1$.
\item $v \in\shortX$ if $v$ is alternating, 
  $\Identity{\comp\opa}$ or $\Identity{\comp\opb}$  appear among
  the
  labels, $\genof{\Identity\opa}$ or $\genof{\Identity\opb}$ do
  not appear as labels, and all of the labels in~$v$
  have \branchcount~$\le 1$.
\end{itemize}
  Every $v\in\longX$ is an alternating composite of elements of the
  form $\genof l$, $l\in\comp\opa\seqcoprod\comp\opb$, such that
  no $l$ is $\Identity{\comp\opa}$ or $\Identity{\comp\opb}$ and at
  least one $l$ is itself a nontrivial composite (in $\comp\opa$ or
  $\comp\opb$). Every element in $\shortX$ is an alternating composite
  of the elements $\genof{\Identity{\comp\opa}}$ and
  $\genof{\Identity{\comp\opb}}$ (at least one of which must appear)
  with elements of the form $\genof{\genof a}$ where
  $a\in\opa\seqcoprod\opb$ and $a$ does not equal $\Identity\opa$ or $\Identity\opb$. 

  It may be helpful in understanding these definitions is to think of
  the case in which \opa/ and \opb/ are graded discrete objects
concentrated at level~1, equivalently,  classic monoids
(\S\ref{CProofOfApplicationTwo}). The free  monoid  on a
set $A$ has as elements the sequences
$(a_1,\ldots, a_k)$ ($k\ge0$) in which the terms $a_i$
are taken from $A$  (in our language, this is a labeling of a linear
tree: the $a_i$'s are the labels, and the sequence length $k$  is
the \branchcount).  Multiplication is given by
concatenation; the empty sequence is allowed, and serves
as the identity element. Then $X_0$ has as elements the
sequences above in which each $a_i$ is taken
either from \opa/ or from \opb/. Similarly $X_1$ is made up of
sequences $v=(u_1,\ldots,u_k)$, ($k\ge0$), such that each $u_i$
itself is either a sequence of elements from \opa/ or a sequence of
elements from \opb. Note that the empty sequence of elements in
\opa/ ($\Identity{\comp\opa}$) is to be distinguished from the empty
sequence of elements in $\opb$ ($\Identity{\comp\opb}$). The map $d_0$
concatenates the sequences $u_i$ together; this usually results in a longer
sequence, but not always, because both varieties
$\Identity{\comp\opa}$ and $\Identity{\comp\opb}$ of empty sequence
collapse when the
concatenation is performed. The map $d_1$ contracts each 
sequence $u_i$ by multiplying it out 
(either in \opa/ or \opb, as appropriate) to a single element $\bar u_i$, 
thus obtaining a sequence $(\bar u_1,\ldots,\bar u_k)$ in $X_0$.
A
sequence $v$ belongs to $\longX$ if the $u_i$ alternate in being
sequences from
\opa/ and sequences from \opb/, no $u_i$ is an empty sequence, and at least one
$u_i$ is a sequence of length $>1$. This guarantees that the length of
$d_0(u)$ is greater than the length of $u$, while the length of
$d_1(u)$ equals the length of $u$. A sequence $u$ belongs to $\shortX$
if the $u_i$ again alternate, at least one $u_i$ is an empty sequence,
each of the $u_i$ has length $\le1$, and no $u_i$ is a length~1
sequence consisting solely of $\Identity\opa$ or $\Identity\opb$. This
guarantees among other things that the length of $d_0v$ is less than the length of $v$,
while the length of $d_1(v)$ equals the length of $v$. These length
(\branchcount) considerations play a role in the proof of \ref{CoprodPushout}.

Given \agradedsubspace/ $Y\subset X_1$, and $l,l'\in X_0$, write $l\simsub Y l'$ and say that $l$ is
congruent to $l'$ mod~$Y$ if
$l$ and $l'$ have the same image in the \nspaces-coequalizer  of $(d_0\vert_Y,
d_1\vert_Y)$.   Say that $Y$ is \emph{saturated} if the equivalence
  relation $\simsub Y$ respects operad composition in~$X_0$.
Let $\bothX=\longX \gradedspacesum \shortX$.
\end{numbered}

\begin{lem}
  The \gradedsubspace/ $\bothX\subset X_1$ is saturated.
\end{lem}

\begin{proof}
  This is a straightforward calculation, but writing it out in detail
  is tedious. Each label in $X_0$ is either
  from \opa/ or from \opb. Congruence mod~\longX/ allows adjacent
  labels of the same color to be composed, so that every simplex of
  $X_0$ is congruent mod~\longX/ to a unique alternating element. Now
  observe that for $v=\genof{\Identity{\comp\opa}}$, for instance,
  $d_0v=\Identity{X_0}$ while $d_1v=\genof{\Identity\opa}$. This implies
  that additionally allowing congruence mod $\shortX$ permits
  leaving out terms of the form $\genof {\Identity\opa}$ or
  $\genof{\Identity\opb}$ from an alternating element of~$X_0$. Of
  course, excising such a term renders the element nonalternating, and
  so it must be reduced further by congruence mod~\longX. (This
  reduction might introduce other terms of the form $\genof
  {\Identity\opa}$ or $\genof{\Identity\opb}$ since some composite of
  elements in $\opa_1$ for instance, might be $\Identity\opa$, but the
reduction process is bound to stop eventually.) The upshot is that
every labeled tree in $X_0$ is congruent mod~$\bothX$ to a unique
alternating labeled tree which contains no
labels of the form $\Identity\opa$ or $\Identity\opb$. Compatibility with
composition is not hard to check inductively.
\end{proof}

Let $\bothd_0$ and $\bothd_1$ denote the restrictions of $d_0$ and
$d_1$ to $\bothX$. In the statement of the following lemma,
coequalizers are to be computed in~\nspaces.

\begin{lem}\label{SmallCoequalizer}
 The natural map from the coequalizer of $(\bothd_0,\bothd_1)$ to the
 coequalizer of $(d_0,d_1)$ is an isomorphism of \gradedspaces.
\end{lem}
\begin{proof}
  Since $\bothX$ is saturated, it is enough to prove that if $v\in X_1$
  has \branchcount~1, then $d_0v$ is congruent to $d_1v$ mod~$\bothX$.
  But this is obvious: either $v=\genof{\Identity{\comp\opa}}$ or
  $\genof{\Identity{\comp\opb}}$ (in which case $v\in \shortX$),
  $v=\genof{\genof a}$ for $a\in\opa\seqcoprod\opb$ (in which case
  $d_0v=d_1v$), or $v=\genof l$ for some
  $l\in\comp\opa\seqcoprod\comp\opb$ of complexity~$>1$, in which case
  $v\in \longX$.
\end{proof}

\begin{titled}{Proof of \ref{CoprodPushout}}
  Let
  $f$ (resp.~$g$) be the map $\bothX\to X_0$  given by $d_0$ on \longX/ and $d_1$
  on~\shortX/ (resp. $d_1$ on \longX/ and $d_0$ on \shortX). By
  \ref{SmallCoequalizer},  $\opa\mycoprod\opb$ is the coequalizer of
  $(f,g)$. It is easy to check that $f$ raises \branchcount/ on
  \longX/ and preserves \branchcount/ on \shortX, whereas $g$
  preserves this invariant on \longX/ and lowers it on
\shortX. The effect of the coequalizer is thus to identify trees of high
\branchcount/ in the image of~$f$ with trees of lower \branchcount/  in
the image of~$g$.   Observe that a simplex of $\Treelabpair\opa\opb{k}$ is
  in $f(\bothX)$  if and only if it belongs to
  $\collTreelab\opa\opb{k}$. The first statement of the proposition now
  follows from relatively routine considerations. The fact that the pushout is a homotopy pushout
  is a consequence of the fact that the upper arrow is a cofibration.
\qed
\end{titled}

\begin{titled}{Proof of \ref{PointedFreeInvariance}}
The fact that $\freeop$ preserves equivalences follows from
\ref{FreeOperadDescription}. To see that $\freeopp$ also preserves equivalences, note that if $X$ is a pointed graded space, there is an exact augmented
split
coequalizer diagram $X_{++}\ReflexArrow X_+\to X$ in the category of
pointed \gradedspaces, obtained by combining the forgetful functor from
pointed objects to unpointed ones with its left adjoint
$(-)_+$. Applying $\freeopp$, using the fact that (as a left adjoint)
$\freeopp$ commutes with colimits, and observing that for, any~$Y$,
$\freeopp(Y_+)=\freeop(Y)$, gives an augmented reflexive coequalizer
diagram
\[
        \freeop(X_+)\ReflexArrow \freeop(X)\to \freeopp(X)
\]
which is exact both in \operads/ and in \nspaces. One of the arrows
$\freeop(X_+)\to\freeop(X)$ in this diagram preserves \branchcount,
the other lowers it. We leave it to the reader to exploit this in
order to give an inductive homotopy invariant construction for
$\freeopp(X)$ along the lines of \ref{CoprodPushout}.
\qed
\end{titled}

\section{Resolutions of bimodules and operads}\label{CResolutions}

\renewcommand{\catM}{\nspace}
\renewcommand{\catMon}{\operads}
In this section we describe a particular operad resolution that
will play a key role in \S\ref{CMapsBetweenOperads}. This resolution
begins as a resolution of a bimodule (\ref{OperadBimodRes}), but it is promoted to a
resolution of an operad by application of an enveloping functor
(\ref{SpecialDistinctionPreserved}).

First, the enveloping functor. 
From now on in this section, \opa/ and
\opb/ are chosen operads. The category \catBiabp/ is the category of pointed
\bimod\opa\opb s; an object $X$ of this category has a basepoint in
level~1, or equivalently, is supplied with a bimodule map
$\opa\ctensor\opb\to X$.  Similarly, \catMonab/ is the category
$\Undercat{(\opa\mycoprod\opb)}{\catMon}$ of operads under $\opa$ and
$\opb$.  Retaining the identity element as
basepoint gives a forgetful functor $\forget\colon\catMonab\to\catBiabp$.
The
following proposition is elementary (it's easy to describe values of
the left adjoint by generators and relations). 

\begin{prop}\label{EnvelopingOperadExists}
  The forgetful functor $\forget\colon\catMonab\to\catBiabp$ has a
  left adjoint $\Envel$.
\end{prop}

From now on in this section, \opa/ is a fixed operad.
The \emph{Hochschild resolution}  of
\opa/ as a bimodule over itself is the simplicial pointed bimodule
$\hres(\opa)$ 
with  \[
    \hres_n(\opa) = \circpower\opa {(n+2)}\,.
  \]
  The object on the right  is a composition power
  of~\opa, with the evident left and right actions.
  The face map $d_i$ is given by using
  operad multiplication to
combine factors $i+1$ and $i+2$ in $\circpower\opa{(n+2)}$, and
  the degeneracy map $s_i$ is given by using the unit inclusion
  $\cunit\to\opa$ to insert the unit
   between factor $i+1$ and factor $i+2$.
The images of
$\hres_0(\opa)$ under the degeneracy maps provide the
necessary basepoints.

The diagonal $\dhres(\opa)$ is the \bimod\opa\opa{} whose
$n$-simplices are the
$n$-simplices of $\hres_n(\opa)$. Operad multiplication maps
$\circpower \opa n\to\opa$ induce a bimodule map
$\daug:\dhres(\opa)\to\opa$, and
the basepoints $\opa\ctensor\opa\to\hres_n(\opa)$ pass to a basepoint 
$\opa\ctensor\opa\to \dhres(\opa)$.

\begin{prop}\label{OperadBimodRes}
  The  bimodule $\dhres(\opa)$ is a
cofibrant object of  \catBiaap, and the map $\dhres(\opa)\to\opa$
above is an
  equivalence. 
\end{prop}

\begin{proof}
  The first statement amounts to a claim that
  $\opa\ctensor\opa\to\dhres(\opa)$ is a cofibration of \bimod\opa\opa
  s, and this is immediate from the description of
  cofibrations in \ref{DescribeCofibrations}. The second one is as
  usual a consequence of the fact that, after forgetting from
  simplicial bimodules to simplicial \gradedspaces, there are
  extra degeneracy maps $s_{-1}:\hres_n(\opa)\to\hres_{n+1}(\opa)$
  which insert the unit as the first composition factor.
\end{proof}

 Applying the functor $\Envel$ (\ref{EnvelopingOperadExists}) 
to $\hres(\opa)$ degree by degree gives a simplicial operad
under $\opa\mycoprod\opa$. The diagonal
$\dehres(\opa)$ is then an operad under $\opa\mycoprod\opa$, and it
too is also supplied with a map $\dehres(\opa)\to\opa$ which factors
the fold map $\opa\mycoprod\opa\to\opa$. This last is a consequence of
the easily verified fact that if $\opa$ is treated as a pointed
bimodule over itself in the natural way, $\Envel(\opa)$ is isomorphic
as an operad under $\opa\mycoprod\opa$ to $\opa$ itself, supplied
with the fold map $\opa\mycoprod\opa\to\opa$.  We now observe that the
same result can be obtained by applying $\Envel$ to the diagonal bimodule
$\dhres(\opa)$. 

\begin{prop}\label{EDefinedSimplicially}
  There two operads
  $\Envel(\dhres(\opa))$ and $\dehres(\opa)$ are isomorphic in a
  natural way (in the category of operads
  under $\opa\mycoprod\opa$ and over $\opa$).
\end{prop}

\begin{proof}
  Let $\opa_k$ be the operad of  sets obtained by
  taking the dimension~$k$ simplices of \opa, and $\Envel_k$ the left
  adjoint to the forgetful functor from operads under
  $\opa_k\mycoprod\opa_k$ to \bimod{\opa_k}{\opa_k}s. It is
  elementary to see that the functor $\Envel$ is cobbled together from
  the functors $\Envel_k$; in other words, if $X$ is a pointed
  bimodule over $\opa$ and $X_k$ is the bimodule over $\opa_k$
  obtained by taking the collection of $k$-simplices in $X$, the
  $\Envel(X)$ is an operad in simplicial sets with
  $\Envel(X)_k=\Envel_k(X_k)$. The result follows immediately.
\end{proof}

\begin{prop}\label{SpecialDistinctionPreserved}
  The operad $\dehres(\opa)$ is a cofibrant object of $\catMonaa$, and
  the map $\dehres(\opa)\to\opa$ is an equivalence.
\end{prop}

\begin{numbered}{Remark} 
  \label{DiagonalPrinciple}\label{SimplicialDot}
  The \emph{diagonal principle} which figures in the following proof
  states that if $X\to Y$ is a map of simplicial (graded) spaces which
  is an equivalence in each simplicial degree, then the induced map
  $\diag X\to\diag Y$ is an equivalence \cite[IV.1.7]{goerssJardine}.

  If  $S$ is a simplicial set and $A$ is an object in some
  category with coproducts, we let $S\simptensor A$ denote the
  coproduct of $S$ copies of $A$. This construction is functorial in
  $S$, and so for a simplicial set $K$ there is a simplicial object
  $K\simptensor A$ which in simplicial degree~$n$ consists of
  $K_n\simptensor A$. If $A$ is itself a simplicial object, then
  $K\simptensor A$ denotes the diagonal of the bisimplicial object
  $\{K_i\simptensor A_j\}$. 
  \end{numbered}

\begin{titled}{Proof of \ref{SpecialDistinctionPreserved}}
  The functor $\Envel$ is left adjoint to a functor which preserves
  fibrations and equivalences, and so
  $\Envel$ preserves cofibrant objects and equivalences between
  cofibrant objects (\ref{QuillenPairs}). The statement
  that $\dehres(\opa)$
  is cofibrant now follows from \ref{EDefinedSimplicially} and
  \ref{OperadBimodRes}.

  Say that the operad \opa/ is \emph{good} if $\dehres(\opa)\to\opa$
  is an equivalence. 
  There  is an isomorphism
  \[
   \Envel\hres_n(\opa) \cong \Envel\circpower\opa {(n+2)} \cong
                      \opa\mycoprod   \freeopp(\circpower\opa
                     n)\mycoprod\opa
   \]
   and so it follows from \ref{CoproductInvariance} and
   \ref{PointedFreeInvariance} that $\Envel\hres_n(\opa)$  depends on~\opa/ in a
   homotopy invariant way. Hence (by the diagonal principle) if $\opa\to\opb$
   is an equivalence, then $\opb$ is good if and only if $\opa$ is. In
   particular, in order to check whether $\opa$ is good, we may assume
   that $\opa$ is cofibrant as an operad (\ref{DescribeCofibrations}).
   Another application of the diagonal principle (cf. the proof of
   \ref{EDefinedSimplicially}) shows that $\opa$ is good if for each
   $n\ge0$ the discrete operad $\opa_n$ is good, where $\opa_n$ is the
   operad of $n$-simplices in~\opa. The upshot is that in order to
   prove that any operad \opa/ is good, it is enough to treat the
   special case in which \opa/
   is the free operad $\freeop X$ on \agradedset/ %
   (discrete \gradedspace)~$X$.

\newcommand{\opad}{\Math{\bar{\opa}}}
\newcommand{\xd}{\Math{\bar X}}

  Given $\opa=\freeop X$, we identify
  $\hres(\opa)$ with $\dhres(\opa)$, since $\hres(\opa)$ is constant in
  the internal simplicial direction. We will construct a cofibrant
  object $J$ of \catBiaap/ together with an equivalence
  $J\to\hres(\opa)$, such that $\Envel(J)$ is clearly equivalent
  to~\opa. The proposition will follow from the fact that $\Envel$
  preserves equivalences between cofibrant bimodules (\ref{AutoQuillenPair}).

  There is a natural map $X\to\opa$ of \gradedspaces/ (the inclusion of generators),
  and so a forgetful functor  $\Psi$ from pointed \bimod\opa\opa s to
  the category of
  \gradedspaces/ under $X\seqcoprod X$. Let $G$ be left adjoint
  to~$\Psi$.  Applying
  $(-)\simptensor X$ to the 
  inclusion $\partial\Delta[1]\to\Delta[1]$ 
  gives $\Delta[1]\simptensor X$ the structure
  of \agradedspace/ under $X\seqcoprod X$. Since 
 $\Delta[1]\simptensor X$ is free as a simplicial \gradedset/ on a
  single copy of the \gradedset/ $X$ in dimension~$1$, there is  a map $\Delta[1]\simptensor
  X\to\Psi\hres(\opa)$  determined by 
  \[   X\to\hres_1(\opa)\cong \opa\ctensor\opa\ctensor\opa, \quad
  x\mapsto \etree (\genof x
  (\etree,\ldots,\etree))\]
   where the number of copies of $\etree$ on the right equals the
  level of~$x$.
   The object $J$ is $G(\Delta[1]\simptensor X)$; adjointness gives a
    map $J\to\hres(\opa)$. We claim that this map is an equivalence,
    that $J$ is cofibrant, and that $\Envel(J)\sim\opa$.

  By inspection, in dimension~$n$  the object $J$ is isomorphic to
  $\opa\ctensor ( X^{\seqcoprod n})_+\ctensor\opa$, where $(-)_+$
  signifies adding a disjoint basepoint at level~1; from this it easily follows
  that $J$ is cofibrant as a pointed bimodule
  (\ref{DescribeCofibrations}). Again by inspection, $\Envel(J)$ is
  isomorphic to $\Delta[1]\simptensor\opa$ (calculated in \operads) or
  equivalently to $\freeop(\Delta[1]\cdot X)$; since $\Delta[1]$ is
  simplicially contractible to either boundary vertex, $\Envel(J)$ is
simplicially contractible to~\opa.

It remains to show that $J\to \hres(\opa)$ is an equivalence. To
do so we prove that the map $J\to\opa$, given by the multiplication
map in degree $0$, is an equivalence. Since this map clearly factors
through the equivalence $\hres(\opa)\to\opa$, we can then
conclude that $J\to \hres(\opa)$ is an equivalence as well.
It is
  easy to check that $\pi_0 J$, which is the coequalizer of the two
  maps $d_0, d_1:\opa\circ X_+\circ\opa\to\opa\circ\opa$ is in fact
  isomorphic to $\opa$ via the operad multiplication
  $\opa\ctensor\opa\to\opa$, so it is enough to check that for each
labeled tree $l\in\opa=\Tof(X)$, the component $C(l)$ of $J$ corresponding
  to $l$ is contractible. We show this by induction on the
  \branchcount/ of~$l$. The complex $C(e)$ is isomorphic to
$\Delta[0]$. Suppose that $l=\treeof(x;l_1,\ldots,l_n)$. The
$k$-simplices of $C(l)$ prescribe ways of splitting~$l$ as a
three-fold composition 
\[l= p_0(\genof{x_1}(p_1^1,\ldots,p_{j_1}^1),\ldots,\genof{x_m}(p_1^m,\ldots,p_{j_m}^m))\,
\]
where $p_0$ and the $p_a^b$'s are elements of \opa, $m$ is the level
of $p_0$, the $x_i$'s belong to $(X^{\seqcoprod k})_+$, and $j_i$ is
the level of $x_i$. (The rogue element $\genof{+}$, where $+$ is
the disjoint basepoint, is to be treated in the
composition formula as the identity element~$e=\Identity{\opa}$.) Consider the following two
possibilities. If $p_0$ is \emph{not}~$e$, then this simplex
lies in the image of the monomorphism
\[ 
      C(l_1)\times\cdots\times C(l_n)\to C(l)\quad\quad
   (\sigma_1,\ldots,\sigma_n)\mapsto \genof x(\sigma_1,\ldots,\sigma_n)\,,
\]
where $\genof x(\sigma_1,\ldots,\sigma_n)$ involves the left action
of $\opa$ on $J$.
If $p_0=e$, then the simplex lies in the image of
\[   
      \Delta[1]\to C(l) \qquad\qquad\quad \quad\Delta_1\mapsto e\circ
       x\circ(l_1,\ldots,l_n)\,,
\]
where $\Delta_1$ is the generating one-simplex of
   $\Delta[1]$. These two images are contractible (the first by induction), they cover $C(l)$,
 and they overlap
in the zero-simplex $\genof x\circ \genof{+}\circ(l_1,\ldots,l_n)$. It
   follows that $C(l)$ is contractible.
   \qed
\end{titled}
\section{Distinguished objects and the proof of \ref{ApplicationOne}}\label{CDistinguished}\label{CMapsBetweenOperads}

\renewcommand{\catM}{\nspace}
\renewcommand{\catMon}{\operads}

In this section we prove a key technical result
(\ref{EnvelhInvariance}) which leads at the end of the section to a
proof of~\ref{ApplicationOne}. We continue to use the notation of
\S\ref{CResolutions}, with the convention that \opa/ and \opb/ are
fixed operads.
An object $X$ of
\catBiabp/ is \emph{distinguished} if the natural map $\opb\to X$ is an
equivalence, and an object \opaa/ of \catMonab/ is distinguished if
$\opb\to \opaa$ is an equivalence.
Note that both \catBiabp/ and \catMonab/ inherit model structures
from~\ref{RezkModelTheorem}. The functor $\Envelh:\catBiabp\to\catMonab$ is the left derived
functor (\ref{AutoQuillenPair}) of $\Envel$ (\ref{EnvelopingOperadExists}) .

\begin{prop}\label{EnvelhInvariance}
  The functor $\Envelh$ preserves distinguished objects.
\end{prop}

We will build up to this in stages. It is convenient to denote a pointed
\bimod\opa\opb{} $X$ by writing out the  triple $T=\triple\opa
X\opb$, where the basepoint $\cunit\to X$ is understood. A morphism of
triples consists of three morphisms, two of operads and one of
pointed bimodules, which are compatible in the obvious sense. The
morphism is an \emph{equivalence} if all of its constituents are
equivalences. The functors $\Envel$ and  $\Envelh$ extend to 
functors on the category of triples.

\begin{prop}
 \label{EquivalencePreservesDistinction}
 If $f:T\to T'$ is an equivalence of triples, then  $\peh(f)$  is an equivalence
 of operads. In particular, if $\peh(T)$ is distinguished, so is $\peh(T')$.
\end{prop}

\begin{proof}
  Write  $T=\triple\opa X\opb$ and
 $T'=\triple{\opa'}{X'}{\opb'}$.
  Let $f^*:\PointedBiModuleCat{\nspaces}{\opa'}{\opb'}\to \PointedBiModuleCat{\nspaces}{\opa}{\opb}$ be
  the restriction functor and $f_*$ its left adjoint. According to \cite[8.6]{rezkProperModel}, the
   pair $(f_*,f^*)$ is a Quillen equivalence. Similarly, let
  $g$ be the map
  $\opa\mycoprod\opb\to\opa'\mycoprod\opb'$,
  $g^*:\BiModuleCat{\operads}{\opa'}{\opb'}\to \BiModuleCat{\operads}{\opa}{\opb}$ the restriction map, and $g_*$ its left adjoint. The map $g$ is an
  equivalence (\ref{CoproductInvariance}) and so according to
  \cite[2.7]{rezkProperModel} it follows from the fact that \operads/
  is left proper (\ref{CoproductInvariance}) that the pair $(g_*,g^*)$
  is a
  Quillen equivalence. It is easy to see that $f_*(X\cof)\to X'$ is an
  equivalence (cf. \cite[pf. of
  8.6]{rezkProperModel}), and so, since $f_*(X\cof)$ is a cofibrant
  object equivalent to~$X'$,  $\Envelh(X')\sim \Envel
  f_*(X\cof)$. Uniqueness of adjoints implies that
  $\Envel f_*\cong g_*\Envel$, and the proposition follows directly.
\end{proof}

A triple $\triple\opa\opb\opb$ is \emph{obtained from a homomorphism
  $f:\opa\to\opb$} if the basepoint in $\opb$ is the operad unit, the
right action of $\opb$ on itself is the usual one, and the left action
of $\opa$ on $\opb$ is obtained by composing the usual left action
of \opb/ on itself with the morphism~$f$. Such a triple is automatically
distinguished.

\begin{prop}\label{MorphismBimoduleDistinguished}
  If a triple $T$ is obtained from
  $f:\opa\to\opb$, then $\peh(T)$ is distinguished.
\end{prop}

\begin{lem}\label{RestrictedCofibration}
  A cofibration $X\to Y$ of cofibrant  \bimod\opa\opb s
  is also a cofibration of cofibrant right \opb~modules. The same
  statement holds for pointed bimodules.
\end{lem}

\begin{proof}
  According to \ref{DescribeCofibrations}, up to retracts 
  the hypotheses amount to
  the condition that there are free degeneracy objects $A\subset
  B\subset Y$ in the category of \gradedsets/ such that $A\subset X$,
  and such that
  the natural maps $\opa\ctensor A\ctensor\opb\to X$ and $\opa\ctensor
  B\ctensor\opb\to Y$ are isomorphisms of degeneracy objects. 
  (Note
  that as in \cite[pf. of 6.2(2)]{rezkProperModel}, it is automatic
  that any basis for $A$ extends to a basis for $B$.)  
In the
  pointed case, the degeneracies of the basepoint are required to lie
  in~$A$. It is enough to show that $\opa\circ A$ and $\opa\circ B$
  are free degeneracy objects in the category of \gradedsets/ such
  that some basis for $\opa\circ A$ can be extended to a basis for
  $\opa\circ B$. In light of the nature of the composition product
  (\ref{DefineCompositionProduct}) and of the fact that as a
  degeneracy diagram each simplicial set $\opa_n$ is free
  \cite[\S6]{rezkProperModel}, the proof comes down to the remark that
  if $S\subset T$ is an inclusion of free degeneracy diagrams of sets,
  and $R$ is another such free degeneracy diagram, then $S\times R\to T\times R$
  is also an inclusion of free degeneracy diagrams. This is proved by
  observing 
  that if $D_i$ is a free degeneracy diagram on a
  element of dimension $i$ then $D_i$ is isomorphic to the diagram of
  simplices in $\Delta[i]$ which do not lie on $\partial\Delta[i]$;
  it follows that  $D_i\times D_j$ is the free degeneracy diagram formed
  by the simplices of $\Delta]i]\times\Delta[j]$ which do not lie on
  the boundary of this product. A basis for $D_i\times
  D_j$ can be described in terms of shuffles \cite[6.5]{maySimplicial}.
\end{proof}

\begin{titled}{Proof of \ref{MorphismBimoduleDistinguished}}
  Let
  $f^*:\PointedBiModuleCat{\nspace}{\opa}{\opb}\to\PointedBiModuleCat
  {\nspaces}{\opa}{\opa}$ be the restriction functor and $f_*$ its
  left adjoint. Similarly, let
  $g^*:\BiModuleCat{\operads}{\opa}{\opb}
  \to\BiModuleCat{\operads}{\opa}{\opa}$ be the restriction functor and
  $g_*$ its left adjoint. As in the proof of
  \ref{EquivalencePreservesDistinction}, we have $\Envel f_*\cong
  g_*\Envel$. The functors $f_*$ and $g_*$ preserve cofibrations and
  equivalences between cofibrant objects (\ref{QuillenPairs}); moreover, in view of the
  formula  $f_*(X)\cong X\ctensor_{\opa}\opb$ \cite[4.4,
  4.7]{rezkProperModel}, $f_*$ is effectively left adjoint to the
  restriction functor from right \opb-modules to right \opa-modules,
  and so preserves equivalences between objects which are cofibrant as
  right \opa-modules. 

   Let $\opa\cof$ be a cofibrant replacement for $\opa$ as
  a pointed \bimod\opa\opa. Since $\opa$ is cofibrant as a right $\opa$-module, the observation above implies that
the map $f_*(\opa\cof)\to
  f_*(\opa)\cong\opb$ is an equivalence, so  $\Envelh(\opb)\sim\Envel
  f_*(\opa\cof)\cong g_*\Envel(\opa\cof)$. 
  Moreover,  $\Envel(\opa\cof)$ is
  distinguished (\ref{EDefinedSimplicially}, 
  \ref{SpecialDistinctionPreserved}).
   Consider the diagram
  \[
  \xymatrix {\opa\ar[r]^{\text{in}_2}\ar[d] &\opa\mycoprod\opa \ar[r]^j\ar[d] &
    \Envel(\opa\cof)\ar[d]\\
             \opb \ar[r] & \opa\mycoprod\opb \ar[r] & g_*\Envel(\opa\cof)
}
  \]
in which the upper composite is an equivalence. 
The small squares are cocartesian as well as homotopy cocartesian \cite[(II.8.14)$\op$~ff.]{goerssJardine}.
The homotopy cocartesian property for the left-hand square is a
consequence of \ref{CoproductInvariance} and for the right-hand square
it is a consequence
of the definition \cite[(II.8.14)$\op$]{goerssJardine} and the fact
that $j$ is a cofibration of operads ($\Envel$ preserves cofibrant
objects). By \cite[(II.8.22)$\op$]{goerssJardine} the large rectangle
is homotopy cocartesian, and the result follows easily.
\qed
\end{titled}

Recall (\ref{GradedTensorProduct})  that category of right \opb-modules is closed under the graded
cartesian product operation~$\gtensor$. The \emph{endomorphism operad}
$\eend_{\opb}(M)$ of a right \opb-module $M$ is
given by
\[
   \eendb(M)\Of n = \Map_{\opb}(M^{\gtensor n}, M)\,.
\]
The subscript $n$ on the left denotes the level of the \gradedspace/
forming the operad, and the object on the right  is a simplicial set of right
\opb-module maps.  The object $\eend_{\opb}(M)$ is an
operad, and $M$
is naturally an
\bimod{\eendb(M)}\opb. 

\begin{rem}\label{EndomorphismsOverSelf}
  In the setting above,  the \bimod\opa\opb{} structures on $M$ which extend the
  given right \opb/ structure are in bijective correspondence with
  operad homomorphisms $\opa\to\eend_{\opb}(M)$.
 If  \opb/ is treated as a right module over itself in the
usual way, then  the left action of \opb/ on itself induces an
isomorphism of operads $\opb\cong\eend_{\opb}(\opb)$. This follows
from the fact that $\opb^{\gtensor n}\cong (*_n)\ctensor\opb$ is the free right
   \opb-module on a single point $(*_n)$ at level~$n$
\end{rem}

If $f:M\to N$ is a map of right \opb-modules, the endomorphism operad
$\eend(f)$ of $f$ is constructed by letting $\eend(f)\Of n$
be the simplicial mapping space of right \opb-module maps $f^{\gtensor
  n}\to f$; this is the space of all commutative diagrams
\[
\xymatrix { M^{\gtensor n} \ar[r]^{f^{\gtensor n}}\ar[d] & N^{\gtensor
    n}\ar[d]\\
     M \ar[r]^f & N}
\]
in which the vertical maps respect the right \opb-actions. There are
natural operad maps $\eend(f)\to\eend(M)$ and $\eend(f)\to\eend(N)$.

 \begin{lem}\label{EndosTheSame}
  If $f:\opb\to N$ is an acyclic cofibration of right
  \opb-modules such that $N$ is fibrant, then
 the natural operad maps $\eend_{\opb}(f)\to\eend_{\opb}(\opb)$ and
  $\eend_{\opb}(f)\to\eend_{\opb}(N)$ are equivalences.
\end{lem}

\begin{proof}
   For each $n\ge0$ there is a fibre square
   \[
     \xymatrix{
          \eend_{\opb}(f)\Of n \ar[r]\ar[d] &  \eend_{\opb}(N)\Of n\cong\Map_{\opb}(N^{\gtensor
            n},N) \ar[d]\\
          \eend_{\opb}(\opb)\Of n\cong\opb_n\ar[r] & \Map_{\opb}(\opb^{\gtensor n},
          N)\cong N_n}
    \]
   For the identifications on the bottom row, recall (\ref{EndomorphismsOverSelf}) that $\opb^{\gtensor n}$ is the free right
   \opb-module on a generator at level~$n$.  The lower map is an
   equivalence, and so it is enough to show that the right vertical
   map is an acyclic fibration, or even that $\opb^{\gtensor n}\to N^{\gtensor n}$ is an acyclic
   cofibration of right \opb-modules.
   This latter map is clearly
   an equivalence. The required cofibration statement follows  from \ref{DescribeCofibrations},
 the distributive formula \ref{DistributiveLaw}, and the fact that the
 product of two free degeneracy diagrams of sets is again a free
 degeneracy diagram (proof of  \ref{RestrictedCofibration}).
\end{proof}

\begin{titled}{Proof of \ref{EnvelhInvariance}}
  Suppose that $M$ is a distinguished \bimod\opa\opb. By the homotopy
  invariance of $\peh$ we can assume that $M$ is fibrant and
  cofibrant as a pointed bimodule, so that in particular the map
  $\opa\circ\opb\to M$ provided by the basepoint is a
  cofibration. Let $f:\opb\to\opa\circ\opb\to M$ be the right
  \opb-module map provided by the basepoint; 
  $f$ is a 
  cofibration of right \opb-modules (\ref{RestrictedCofibration}) and
  the natural operad maps $\eend_{\opb}(f)\to\eend_{\opb}(M)$ and
  $\eend_{\opb}(f)\to\opb$ (\ref{EndomorphismsOverSelf}) are equivalences~(\ref{EndosTheSame}).
  The left action
  of $\opa$ on $M$ gives an operad homomorphism $\opa\to \eend_{\opb}(M)$ which
  we factor as the composite of an acyclic cofibration $\opa\to\opa'$
  and a fibration $\opa'\to\eend_{\opb}(M)$. Let $\opa''$ be given by the
  pullback diagram
  \[
     \xymatrix{\opa''\ar[r]\ar[d] & \opa'\ar[d]\\
                 \eend_{\opb}(f) \ar[r] &\eend_{\opb}(M) }
  \]
  The upper map here is an equivalence because the lower one is.
  There are 
  equivalences of pointed bimodule triples
  \[(\opa,M,\opb)\rightarrow (\opa',M,\opb)\leftarrow
  (\opa'',M,\opb)\leftarrow(\opa'',\opb,\opb)\]
where the action of $\opa''$ on $\opb$ arises from the operad map  $\opa''\to\eend_{\opb}(f)\to\opb$. Since
  $\peh(\opa'',\opb,\opb)$ is distinguished
  (\ref{MorphismBimoduleDistinguished}), so is $\peh(\opa,M,\opb)$
(\ref{EquivalencePreservesDistinction}).
  \qed
\end{titled}

\begin{titled}{Proof of \ref{ApplicationOne}} We have to check Axioms
  I-VI from \S\ref{MonoidalCategories} in the case in which $\origcatM$
  is the category of \gradedspaces/ and $\gentensor$ is the
  composition product. Theorem \ref{LoopedMonoidMapFibration} then yields
  \ref{ApplicationOne}; the assumption that $\oo\Of 1$ is contractible guarantees
  that $\Map_{\nspaces}(\cunit,\oo)\sim\oo\Of 1$ is contractible.
  
  Axiom~I is \ref{RezkModelTheorem}; Axioms~II and~III are
  \ref{EnvelopingOperadExists} and \ref{EnvelhInvariance}
  respectively. Axiom~IV(2) is \ref{CoproductInvariance} and
  Axiom~V(2) follows as in \ref{GeneratorCofibrant} from the fact that
  $\cunit$~is cofibrant as \agradedspace. Axiom~VI(1) is trivial,
  since $\cunit\cof\cong\cunit$. Finally, Axiom~VI(2) is a consequence
  of the fact 
that if \opa/
  is an operad, any cofibrant replacement $\opa\cof$ for $\opa$ as a
  \bimod\opa\opa{} is also cofibrant as a right \opa-module
  (\ref{RestrictedCofibration}): given an operad map $\opa\to\opb$,
  the functor $(\text{--})\ctensor_{\opa}\opb$ is left adjoint to the
  forgetful  functor from right $\opb$-modules to right
  $\opa$-modules; since the forgetful functor preserves equivalences
  and fibration, $(\text{--})\ctensor_{\opa}\opb$ preserves
  equivalences between cofibrant right \opb-modules.
\qed
\end{titled}

\section{\Classicmonoids/ and the proof of \ref{ApplicationTwo}}
\label{CProofOfApplicationTwo}

A \classicmonoid/ is a monoid object in $(\spaces,\times,*)$, i.e.,
a simplicial semigroup with identity, or equivalently, an operad
concentrated at level~1. This last observation allows all of the
results of the previous sections to be applied to
\classicmonoids. With only notational changes, the same results (with
the same proofs) hold for \gradedclassicmonoids/, i.e., monoid
objects in $(\nspaces,\gtensor,\gunit)$. 

The variant of  Theorem
\ref{ApplicationOne} in which \assoc/ is replaced by an arbitrary
operad (the variant actually proved in \S\ref{CMapsBetweenOperads}) now translates to the following statement. Let 
\cgm/ denotes the category of \gradedclassicmonoids. 

\begin{thm}\label{TranslationOfTwo}
  Suppose that $\alpha:G\to H$ is a map of \classicgradedmonoids, such that
  $H\Of0\sim*$. Then there is an equivalence
  \[
     \Omega\Maph_{\cgm}(G,H)_\alpha \sim \Maph_{G-G}(G,H)\,,
  \]
  where on the right $H$ is treated as a \bimod GG via~$\alpha$.
\end{thm}

\begin{rem}
The grading on these objects 
 allows the assumption in \ref{ApplicationOne} that
the target operad is contractible at level~1 (which for ungraded
classic monoids
would 
translate into the unfortunate assumption that the target monoid is contractible)
to be replaced by the 
assumption that the grade~0 constituent of the target graded monoid, in other words
the
constituent  containing the identity element, is contractible.
\end{rem}

A \gradedclassicmonoid/ is simply a left module over the associative operad
\assoc, so specializing
\ref{TranslationOfTwo} to the case $G=\assoc$ gives the following
result. For this statement, \gassocn/  refers to the operad \assoc, treated
as a left module over itself, i.e., as a graded classic monoid.

\begin{thm}\label{SpecialTranslation}
  Suppose that $\alpha:\gassocn\to X$ is a map of  left $\assoc$-modules, and
  that $X\Of0\sim *$. Then there is an equivalence
  \[
    \Omega\Maph_{\assoc}(\gassocn,
    X)_\alpha\sim\Maph_{\gassocn-\gassocn}(\gassocn, X)\,.
  \]
\end{thm}

  The mapping space on the right above is a derived mapping space of \gtensor-bimodules over
  the \gradedclassicmonoid/ \gassocn, where $X$ is treated as an
  \bimod{\gassocn}{\gassocn}{} via~$\alpha$. 
Theorem \ref{SpecialTranslation} is remarkably similar to \ref{ApplicationTwo}; the only
difference being that, in \ref{ApplicationTwo}, $\gassocn$ and $X$ have
additional right $\assoc$-module structures, and the mapping spaces
respect these additional structures. We will proceed to prove
\ref{ApplicationTwo} by showing that introducing  right \assoc-module
structures does not materially change the arguments.  

Some notation will be useful. As usual, \nspaces/ is the category of
\gradedspaces/, and for the purposes below $\RtAmod$ will denote the category of right
\assoc-modules. We use \LeftAmod/ to
denote the
category of left \assoc-modules and
\LeftAAmod/ to denote the category of left \assoc-modules in \RtAmod/ (i.e., the
category of \bimod\assoc\assoc s). An object
of \LeftAmod/ is a graded classic monoid, and an object of \LeftAAmod/
is a graded classic monoid with a right \assoc-action 
compatible with the monoid multiplication. 

If \gcma/ and \gcmb/ belong
to \LeftAAmod, we use the same letters to denote the underlying
objects of \LeftAmod. There are then categories
$\ModuleCat{\nspaces}{\gcmb}$, $\ModuleCat{\RtAmod}{\gcmb}$,
$\BiModuleCat{\nspaces}\gcma\gcmb$,
$\BiModuleCat{\RtAmod}{\gcma}\gcmb$ of right \gcmb-modules or
\bimod\gcma\gcmb s,  as well as pointed
variants $\PointedBiModuleCat{\nspaces}\gcma\gcmb$ and
$\PointedBiModuleCat{\RtAmod}{\gcma}\gcmb$. Keep in mind that these
 latter module
structures refer to the graded cartesian product structure~\gtensor, and that the
basepoint for an object $X$ of
$\PointedBiModuleCat{\RtAmod}{\gcma}\gcmb$, for instance,
is a map $\gunit\to X$, and so amounts to a point at level \emph{zero}
in the
underlying \gradedspace.
Disregarding the right \assoc-module structure gives rise to
forgetful functors: $\RtAmod\to\nspaces$, $\LeftAAmod\to\LeftAmod$, 
$\ModuleCat{\RtAmod}{\gcmb}\to\ModuleCat{\nspaces}{\gcmb}$,
$\BiModuleCat{\RtAmod}{\gcma}\gcmb\to\BiModuleCat{\nspaces}\gcma\gcmb$,
and
$\PointedBiModuleCat{\RtAmod}{\gcma}\gcmb\to\PointedBiModuleCat{\nspaces}\gcma\gcmb$.

\begin{lem}\label{AllColimitsMatch}
All of the above forgetful functors are left adjoints, and so preserve
colimits. All of the functors preserve (and reflect) equivalences, and
preserve cofibrant objects.
\end{lem}

\begin{proof}
  We will show that the functors preserve colimits, and the fact that
  they are left adjoints will be a consequence of the Adjoint Functor
  Theorem \cite[V.6]{maclaneCategories}.  Let $U:\Cat A\to\Cat B$ be
  one of these forgetful functors.  The functor $U$
   preserves colimits if and only if it preserves coproducts and
   reflexive coequalizers; $U$ clearly preserves reflexive
   coequalizers, because as in \ref{ReflexiveCoequalizers} these can
   be computed in~\nspaces.  Say that an object of $\Cat A$ is
  \emph{free} if it is in the image of the left adjoint $\Phi$  to the 
   forgetful functor $F: \Cat A\to \nspaces$. Since every object $X$ of $\Cat A$
   is given  a reflexive coequalizer $(\Phi
   F)^2(X)\ReflexArrow\Phi F(X)$, in order to prove
   that $U$ preserves coproducts it is enough to show that $U$ preserves
   coproducts of free objects. But in each case this is clear by
   inspection, e.g., for $U:\RtAmod\to\nspaces$ and for free objects $X\ctensor \assoc$ and
   $Y\ctensor\assoc$ of \RtAmod,
  \[\begin{aligned}
U(X\ctensor
  \assoc\coprodOver{\RtAmod}(Y\ctensor\assoc))&\cong
  U((X\coprodOver{\nspaces} Y)\ctensor\assoc)\\
  &\cong 
  X\ctensor\assoc\coprodOver{\nspaces}Y\ctensor\assoc\,.
  \end{aligned}
\]
To put these isomorphisms in context, it is useful to observe that the
functor $-\ctensor\assoc:\nspace\to\nspace$ preserves colimits, and that seen as a functor $\nspace\to \nspace_{\assoc}$, it creates colimits as well. 
 It is clear that a map $f$ in $\Cat A$ is an equivalence
 if and only if $U(f)$~is an equivalence. The statement about
 preservation of cofibrant objects follows from \ref{DescribeCofibrations}
 and an argument as above depending on the fact that $U$ preserves
 free objects.
\end{proof}

\begin{titled}{Proof of \ref{ApplicationTwo}}
We have to check
Axioms I-VI from \S\ref{MonoidalCategories} in the case in which
\origcatM/ is the monoidal category $(\RtAmod,\gtensor,\gunit)$; as
above, we know that the axioms hold in $(\nspaces,\gtensor, \gunit)$.
Theorem \ref{LoopedMonoidMapFibration} will then yield
\ref{ApplicationTwo}; note that the assumption in \ref{ApplicationTwo}
that $X\Of0$ is contractible guarantees that
$\Map_{\nspaces}(\gunit,X)$ is contractible. 

Axiom~I
follows as in \ref{RezkModelTheorem} from \cite[7.1]{rezkProperModel},
since all of the categories involved are categories of algebras over
$\mathbb{Z}^+$-sorted theories. The existence of the left adjoint
required in Axiom~III is standard.  Suppose that \gcma/ and \gcmb/
belong to \LeftAAmod/; denote the coproduct in
\LeftAAmod/ of these objects by $\gcma\mycoprod\gcmb$. By \ref{AllColimitsMatch},
there is no harm in using the same notation for the coproduct of the
images of \gcma/ and \gcmb/ under the forgetful functor
$\LeftAAmod\to\LeftAmod$. There is a diagram of functors
\[
\xymatrix{
  \PointedBiModuleCat{\RtAmod}{\gcma}\gcmb
  \ar[d]^{\Envel_1}\ar[r]^{F_1}& {\PointedBiModuleCat{\nspaces}\gcma\gcmb}
  \ar[d]^{\Envel_2}\\
   {\Undercat{\gcma\mycoprod\gcmb}{\LeftAAmod}} \ar[r]^{F_2} &
   {\Undercat{\gcma\mycoprod\gcmb}{\LeftAmod}}
 }
\]
 in which the horizontal arrows are forgetful functors and the
 vertical arrows are the appropriate left adjoints to forgetful functors. Adjointness
 gives a natural transformation  $\Envel_2F_1\to F_2\Envel_1$. Since
 all of these functors commute with colimits (\ref{AllColimitsMatch}),
 to show that the diagram commutes up to natural isomorphism it is
 enough to show that the indicated natural transformation is an
 isomorphism when applied to a free object $\gcma \gtensor (X_+\ctensor\assoc)\gtensor
 \gcmb$ of $\PointedBiModuleCat{\RtAmod}{\gcma}\gcmb$. (See the proof
 of \ref{AllColimitsMatch}; here $X\in\nspaces$ and $X_+$ is obtained
 by adding a disjoint basepoint at level~0 to~$X$.) But the
 compositions $F_2\Envel_1$ and $\Envel_2F_1$ 
 take  this free object  to $\gcma\coprodOver{\Cat
   D}\gcmb\coprodOver{\Cat D}
 (\assoc\ctensor X\ctensor\assoc)$, where $\Cat D$ is respectively
 $\LeftAAmod$ or \LeftAmod;  by \ref{AllColimitsMatch}, the
 distinction between the two types of coproduct is irrelevant. 
 The functor $F_1$ preserves equivalences 
 distinguished objects, and cofibrant objects (\ref{AllColimitsMatch}),
 and so $\Envelh_1$ preserves distinguished objects because
 $\Envelh_2$ does. This verifies Axiom~III. 

 The category \LeftAAmod/ is left proper (and  coproducts in \LeftAAmod/
 preserve equivalences) because of \ref{AllColimitsMatch} and the fact
 that the corresponding statements hold in \LeftAmod. This is
 Axiom~IV(2). Axiom~V(1) is immediate, because $\gunit$ is
 cofibrant as \agradedspace. Axiom~VI(1) is trivial, again because
 $\gunit$ is cofibrant. Finally, Axiom~VI(2) can be derived from
 combining the fact that the appropriate forgetful functor preserves
 coequalizers (\ref{AllColimitsMatch}) with the fact that the
 corresponding axiom holds in \LeftAmod.
\qed
\end{titled}

\section{The main theorem}\label{CFinalProof}

Here we prove Theorem~\ref{MainTheorem}; we continue to use 
notation from \S\ref{CProofOfApplicationTwo}.  As in
\ref{ApplicationTwo}, \gassoc/ denotes \assoc/ in its role as a monoid
in $(\RtAmod,\gtensor,\gunit)$. Suppose that \oo/ is a
reduced operad and that $\om:\assoc\to\oo$ is an operad map. 
Combining \ref{ApplicationOne} with \ref{ApplicationTwo} gives an
equivalence
\begin{equation}\label{DoubleLoopEqv}
   \Omega^2\Maph_{\operads}(\assoc,\oo)_{\om}\sim\maptrih\gassoc\gassoc\assoc(\gassoc,\oo)\,,
\end{equation}
where $\oo$ is treated on the right as an object of
$\BiModuleCat{\RtAmod}{\gassoc}\gassoc$. We will calculate the mapping
space on the right by using a specific resolution of
\gassoc/ in \FunnyBiMods.

The \emph{Hochschild resolution} of $\gassoc$ as a bimodule over
itself in \RtAmod/ is the simplicial object $\hres(\gassoc)$ with
\[
   \hres_n(\gassoc) = \gassoc^{\gtensor (n+2)}\,.
\]
The object on the right is a graded cartesian power of \gassoc, with
the evident right and left actions of~\gassoc; the simplicial
operators $\hres_n(\gassoc)\to\hres_{n\pm1}(\gassoc)$ are as in
\S\ref{CResolutions}. Forgetting the right \assoc-module
structures and the grading of \gassoc/ renders $\hres(\gassoc)$ identical to the resolution of \gassoc/ that
would be obtained from \S\ref{CResolutions} by treating \gassoc/ as an
operad concentrated at level~1 (a classic monoid). Let
$\diag\hres(\gassoc)$ be the diagonal of this simplicial object. There
is a map $\aug:\diag\hres(\gassoc)\to\gassoc$ of objects in
$\FunnyBiMods$, and as in \S\ref{CResolutions}, this map is an
equivalence.

The proof of \ref{MainTheorem} depends on two lemmas.

\begin{lem}\label{MapCosimplicialResolution}
  Let $\cd=\FunnyBiMods$, let $X$ be a simplicial object in \cd, and let
  $Y$  be object of $\cd$. Then there is a natural weak homotopy
  equivalence of spaces
  \[
        \Maph_{\cd}(\diag X, Y) \sim \holim \Maph_{\cd}(X, Y)\,.
  \] 
\end{lem}

The right hand side above is the homotopy limit of the cosimplicial space
obtained by applying $\Maph_{\cd}({\text{--}},Y)$ to $X$ degree by degree. 

\begin{lem}\label{IdentifyCoSimplicial}
  Let $\cd=\FunnyBiMods$. The cosimplicial space
  $\Maph_{\cd}(\hres(\gassoc), \oo)$ is equivalent
  to the cosimplicial space \cso/ associated to \om/ by McClure and
  Smith (\S\ref{CIntroduction}).
\end{lem}

\begin{titled}{Proof of \ref{MainTheorem}}
Let $\cd=\FunnyBiMods$, and 
let $X=\hres(\gassoc)$ in \cd. Since
$\diag(X)\sim\gassoc$, \ref{DoubleLoopEqv} gives an equivalence
\[
 \Omega^2\Maph_{\operads}(\assoc,\oo)_{\om}\sim\Maph_{\cd}(\diag(X),\oo)\,,
\]
By \ref{MapCosimplicialResolution}, the space on the right is
equivalent to the homotopy limit of a cosimplicial space obtained by
applying $\Maph_{\cd}(\text{--},\oo)$ to $X$ degree by degree. By
\ref{IdentifyCoSimplicial}, this cosimplicial space is equivalent to the one
appearing in \ref{MainTheorem}. \qed
\end{titled}

\begin{titled}{Proof of \ref{MapCosimplicialResolution}} 
According to \cite[18.1.10]{hirschhorn} or more generally
\cite[19.4.4(1)]{hirschhorn}, it is enough to show that $\diag X$
is naturally weakly equivalent to $\hocolim X$. Since \cd/ is a
simplicial model category, this homotopy colimit can be obtained by the formula
of \cite[18.1.2]{hirschhorn}. Of course the constructions which
figure in this formula should be made in \cd, but since the forgetful
functor $\cd\to\nspaces$ respects colimits, the constructions
could equally well be made in~\nspaces. The lemma now follows easily
from the fact that the homotopy colimit of a simplicial object in
\nspaces/ is naturally weakly homotopy equivalent to its diagonal (cf.
\cite[16.11.6 and 18.7.5]{hirschhorn}).
\qed
\end{titled}

\begin{titled}{Proof of \ref{IdentifyCoSimplicial}}
The \gradedspace/ $\hres_n(\gassoc)$ is the
free object of $\BiModuleCat{\nspaces}{\gassoc}{\gassoc}$  on the right \assoc-module
$\gassoc^{\gtensor n}$. Let $*_n$ be the \gradedspace/ which is empty
except for a single point at level~$n$; it is immediate that 
\[ (*_n)\ctensor\assoc\cong\assoc^{\gtensor n}\cong\gassoc^{\gtensor
  n}
\]
is the free right \assoc-module on $*_n$. It follows that
$\hres_n(\gassoc)$ is the free object of \cd/ on the \gradedspace/ $*_n$, so that
$\hres_n(\gassoc)$ is a cofibrant object of \cd/ and
\[
\Maph_{\cd}(\hres_n(\gassoc,\oo)\sim\Maph_{\nspaces}(*_n,\oo)\cong
\oo\Of n\,.
\]
The cosimplicial operators in $\Maph_{\cd}(\hres_n(\gassoc,\oo)$ can
now be identified by inspecting the simplicial operators in the
Hochschild resolution~$\hres(\gassoc)$.
\qed
\end{titled}

\providecommand{\bysame}{\leavevmode\hbox to3em{\hrulefill}\thinspace}

\end{document}